\documentclass[12pt]{article}%
\usepackage{amsmath,amssymb,graphicx,mathrsfs,harvard,fancyhdr,pictex}%
\renewcommand{\mathcal}{\mathscr}
\newcommand{\proof}{{\it Proof.\ }}
\newcommand{\Var}{{\rm Var}}
\newcommand{\Corr}{{\rm Corr}}

\def\tfrac#1#2{{\textstyle\frac {#1}{#2}}}
\newtheorem{theorem}{Theorem}[section]
\newtheorem{rema}[theorem]{Remark}

\newenvironment{rem}{%
\begin{rema}\bf \rm }{\end{rema}
}

\newtheorem{definition}[theorem]{Definition}

\newtheorem{lemma}[theorem]{Lemma}

\newtheorem{proposition}[theorem]{Proposition}

\begin{document}
\date{}
\title{HIGH FREQUENCY ASYMPTOTICS FOR WAVELET-BASED TESTS FOR GAUSSIANITY AND
ISOTROPY ON THE TORUS }
\author{Paolo Baldi\\
{\it Dipartimento di Matematica,
Universit\`a di Roma {\sl Tor Vergata}, Italy}\\
G\'erard Kerkyacharian\\
{\it Laboratoire de Probabilit\'es et
Mod\`eles Al\'eatoires}\\
{\it and Laboratoire MODALX, Paris, France}\\
Domenico Marinucci\\
{\it Dipartimento di Matematica,
Universit\`a di Roma {\sl Tor Vergata}, Italy}\\
Dominique Picard\\
{\it Laboratoire de Probabilit\'es et Mod\`eles Al\'eatoires,
Paris, France}} \maketitle

\noindent{\it Running title} High Frequency Asymptotics on the Torus.
\medskip

\medskip



\begin{abstract}
\noindent We prove a CLT for skewness and kurtosis of the wavelets
coefficients
of a stationary field on the torus. The results are in the
framework of the fixed-domain asymptotics, i.e. we refer to
observations of a single field which is sampled at higher and
higher frequencies. We consider also studentized statistics for the case of
an unknown correlation structure%
. The results are motivated by the analysis of cosmological data
or high-frequency financial data sets, with a particular interest
towards testing for Gaussianity and isotropy.

\end{abstract}

\section{Introduction}\label{intro}
In recent years, a growing interest has been drawn by infill (or
fixed-domain, or high-resolution)  asymptotics,
i.e., the investigation of the limiting behaviour of statistics
of a single observation of a stochastic processes on a fixed time
span or of a random field on a compact space, observed
at a greater and greater resolution as the sample size increases
(see \cite{MR1682572} e.g.).

Such interest was mainly stimulated by various fields of
applications. For stochastic processes defined on subsets of
$\mathbb{R},$ the leading motivations have come from finance,
where it is now customary to observe high-frequency or ultra-high
frequency data sets collected in a short amount of time (even a
single day). In such cases, clearly, the standard asymptotic
framework, envisaging an ever-increasing time span of
observations, can be highly misleading and a fixed-domain approach
seems a valid alternative. Many important contributions have
focussed on asymptotic statistical inference on discretely
observed diffusions, see for instance Kessler and S{\o}rensen
\cite{MR1681700}, S{\o}rensen and Uchida \cite{MR2046817}.

Our interest here, however, is not on diffusions, but rather
on stationary processes, as for instance in Stein \cite{MR1697409}
and \cite{MR1892665}. In this area, strong interest has been
brought in by the cosmological and astrophysical literature, a
particularly active field being the analysis of Cosmic Microwave
Background radiation. The latter can be viewed as a relic of the
Big Bang which provides a snapshot of the Universe some 13,7
billion years ago; as such it is considered a goldmine of
information on fundamental physics. A widely debated issue relates
to the probability law of such radiation: this has fuelled
considerable activity on testing for Gaussianity on isotropic
random fields defined on
$\mathbb{S}^{2}%
$, the sphere in three-dimensional Euclidean space. Wavelet analysis
has also been proposed here, for instance by Vielva et al.
\cite{Vielva}, Jin et al. \cite{Jin}, Cabella et al.
\cite{cabella:063007}. See also \cite{mari2006a} for rigorous
results on non Gaussianity testing, but based on the angular
bispectrum.

Wavelet analysis has proved to be a powerful tool in this kind of
problems, owing to their good localization properties.
To the best of our knowledge, however, no theoretical analysis exists so far
in the literature on the asymptotic (high-resolution) behaviour of the random
wavelets coefficients for random fields on bounded domains.
In this paper, we provide some preliminary
results in this area. In particular, we focus on isotropic random
fields on
the torus $\mathbb{S}^{1};$ we define a suitable wavelet expansion
and we
derive the correlation structure of the random wavelets coefficients. We use
these results to establish a central limit theorem for the skewness and
kurtosis
statistics. These are useful to investigate the Gaussianity of the
field. As explained
before, the asymptotic theory is clearly developed in the high-resolution
sense. The results are then extended to the case where the correlation
structure of the field is unknown and estimated from the data.

Our choice of $\mathbb{S}^{1}$ as the domain of the field is
a
first-step investigation: it is certainly desirable
to extend our results, for instance, to bounded subsets of
$\mathbb{R}$ or to $\mathbb{S}^{2}$; such extensions are currently
under investigation. However also the present results are of
a practical interest. Indeed, random fields on the circle are the
natural environment for many geophysical models, for instance
concerning atmospheric data \cite{NT95}
. Also, it is
not unusual in the CMB\ literature to approximate $\mathbb{S}^{2}$
as the union of copies of $\mathbb{S}^{1}$, the so-called ring-torus approach
\cite{ring-torus}. Finally, our assumptions can be used for models of stationary
time series on a fixed domain.

The plan of the paper is as follows. In \S \ref{petru} we review
the Petrushev construction of needlets on general spaces; in \S
\ref{assu} and \ref{asym} we introduce random fields on the torus
and the associated random wavelets coefficients; a fundamental
bound is established on the covariances of the latter. \S
\ref{tcl} is devoted to a Central Limit Theorem for the Skewness
and Kurtosis statistics of these random wavelets coefficients.
This result is extended to studentized statistics in \S
\ref{stud}, whereas \S \ref{alias} is devoted to the investigation
of the aliasing effect, that is, the discretized sampling of the
continuous random field of interest. \S \ref{montec} presents some
Monte Carlo evidence.
In the
sequel, we use $c$ to denote a positive constant, which need not be the same
from line to line.

\section{Petrushev  and coauthors construction of Needlets}\label{petru}

Frames were introduced in the 1950's \cite{dufschaef}
to represent functions via over-complete sets. Frames including tight frames
arise naturally in wavelet analysis on $\mathbb{R}^{d}$. Tight frames which
are very close to orthonormal bases are especially useful in signal and image processing.

We shall see that the following construction has the advantage of
being easily computable and of producing well localized tight
frames constructed on a specified orthonormal basis.
%
%
\begin{definition}
Let ${\mathbb{H}}$ be a Hilbert space, and $(e_{n})$ a sequence in
${\mathbb{H}};$ $(e_{n})$ is a tight frame (with constant 1) if :
\[
\forall f\in{\mathbb{H}},\;\Vert f\Vert^{2}=\sum_{n}|\langle f,e_{n}%
\rangle|^{2}%
\]
\end{definition}
Let now $\mathcal{Y}$ be a metric space endowed with a finite
measure $\mu$. Assume that the following decomposition holds
\[
\mathbb{L}_{2}(\mathcal{Y},\mu)=\bigoplus_{l=0}^{\infty}H_{l}%
\]
where the $H_{l}$'s are finite dimensional spaces. For the sake of
simplicity, we suppose that $H_{0}$ is reduced to the constants.
Let $L_{l}$ be the orthogonal projection on $H_{l}$ :
\[
\forall f\in\mathbb{L}_{2}(\mathcal{Y},\mu),\qquad L_{l}(f)(x)=\int
_{\mathcal{Y}}f(y)L_{l}(x,y)d\mu(y)
\]
where
\[
L_{l}(x,y)=\sum_{i=1}^{m_{l}}e_{i}^{l}(x)\bar{e}_{i}^{l}(y)
\]
$m_{l}$ is the dimension of $H_{l}$ and $(e_{i}^{l})_{i=1,\dots,m_{l}}$ an
orthonormal basis of $H_{l}$. Let us observe that we have the following
property of the projection operators:
\begin{equation}
\int L_{l}(x,y)L_{m}(y,z)d\mu(z)=\delta_{l,m}L_{l}(x,z)\ . \label{auto}%
\end{equation}
The following construction, also inspired by
\cite{frazjawe}, is based on two fundamental steps: Littlewood-Paley
decomposition and discretization, which are summarized in the two following
subsections.
\subsection{Littlewood -Paley decomposition}
Let $\phi$ be a $C^{\infty}$ function supported in $|\xi|\leq1,$ such that
$1\geq\phi(\xi)\geq0$ and $\phi(\xi)=1$ if $|\xi|\leq\frac{1}{2}$. Let us
define:
\[
a^{2}(\xi)=\phi(\tfrac \xi2)-\phi(\xi)\geq0
\]
so that
\begin{equation}
\forall|\xi|\geq1,\;\qquad\sum_{j}a^{2}(\tfrac \xi{2^j})=1\ . \label{1}%
\end{equation}
Actually in the previous sum all middle terms cancel
telescopically. Let us define the operator
\[
\Lambda_{j}=\sum_{l\geq0}a^{2}(\tfrac l{2^j})L_{l}%
\]
and the associated kernel
\[
\Lambda_{j}(x,y)=\sum_{l\geq0}a^2(\tfrac l{2^j})L_{l}(x,y)=\sum_{2^{j-1}%
<l<2^{j+1}}a^2(\tfrac l{2^j})L_{l}(x,y).
\]
Then it holds:
\begin{proposition}\label{Mdef}%
\begin{equation}
\forall f\in{\mathbb{H}},\text{\ \ \ }~f=\lim_{J\rightarrow\infty}%
L_{0}(f)+\sum_{j=0}^{J}\Lambda_{j}(f) \label{rep}%
\end{equation}%
and, if $M_{j}(x,y)=\sum_{l}a(\tfrac l{2^j})L_{l}(x,y)$,
\begin{equation}
\Lambda_{j}(x,y)=\int
M_{j}(x,z)M_{j}(z,y)d\mu(z)\ . \label{sqrt}%
\end{equation}
\end{proposition}
\proof%
\begin{equation}
\label{1}L_{0}(f) + \sum_{j=0}^{J} \Lambda_j(f) = L_{0}(f) +
\sum
_{j=0}^{J} \sum_{l} a^{2}(\tfrac  l{2^j})L_{l}(f)  =
\sum_{l}\phi(\tfrac l{2^{J+1}})L_{l}(f).
\end{equation}
Therefore, as $\phi(\tfrac l{2^{J+1}})=1$ as soon as $2^J\ge l$,
\begin{align*}
\big\Vert\sum_{l}\phi(\tfrac l{2^{J+1}})L_{l}(f)-f\big\Vert^{2}  &  =\sum_{l\geq2^{J+1}}\Vert
L_{l}(f)\Vert^{2}+\sum_{2^{J}\leq l<2^{J+1}}\big\Vert L_{l}(f)(1-\phi(\tfrac l{2^{J+1}})\big\Vert^{2}\\
&  \leq\sum_{l\geq2^{J}}\Vert L_{l}(f)\Vert^{2}\longrightarrow0,
\text{ as }J\rightarrow\infty.
\end{align*}
(\ref{sqrt}) is a simple consequence of (\ref{auto}).

\hfill$\square$
\subsection{Discretization}
Let us define
\[
\mathcal{K}_{l}=\bigoplus_{m=0}^{l}H_{m},
\]
and let us assume that some additional assumptions are true:

1. $
f\in\mathcal{K}_{l}\Rightarrow\bar{f}\in\mathcal{K}_{l}%
$

2. $
f\in\mathcal{K}_{l},~~g\in\mathcal{K}_{l}\Rightarrow fg\in\mathcal{K}%
_{l+l}%
$

3. Quadrature formula: for every $l\in\mathbb{N}$ there exists
 a finite subset $\mathcal{X}_{l}\subset\mathcal{Y}$ and positive
 real numbers $\lambda_{\eta}>0$, $\eta\in\mathcal{X}_{l}$, such that
\begin{equation}
\forall f\in\mathcal{K}_{l},~~\int fd\mu=\sum_{\eta\in\mathcal{X}_{l}}%
\lambda_{\eta}f(\eta). \label{quadr}%
\end{equation}
Then the operator $M_{j}$ defined in Proposition \ref{Mdef} is such
that $M_{j}(x,z)=\overline{M_{j}(z,x)}$ and
\[
z\mapsto M_{j}(x,z)\in\mathcal{K}_{2^{j+1}-1}\text{ ,}%
\]
so that
\[
z\mapsto M_{j}(x,z)M_{j}(z,y)\in\mathcal{K}_{2^{j+2}-2}\text{ },
\]
and we can write:
\[
\Lambda_{j}(x,y)=\int M_{j}(x,z)M_{j}(z,y)\,d\mu(z)=\sum_{\eta\in\mathcal{X}%
_{2^{j+2}-2}}\lambda_{\eta}M_{j}(x,\eta)M_{j}(\eta,y)\text{ }.
\]
This implies:
\begin{align*}
\Lambda_{j}f(x)  &  =\int\Lambda_{j}(x,y)f(y)\,d\mu(y)=\int\sum_{\eta
\in\mathcal{X}_{2^{j+2}-2}}\lambda_{\eta}M_{j}(x,\eta)M_{j}(\eta
,y)f(y)\,d\mu(y)\\
&  =\sum_{\eta\in\mathcal{X}_{2^{j+2}-2}}\sqrt{\lambda_{\eta}}M_{j}%
(x,\eta)\int\overline{\sqrt{\lambda_{\eta}}M_{j}(y,\eta)}f(y)\,d\mu(y).
\end{align*}
This can be summarized in the following way, if we set
\[
\mathcal{X}_{2^{j+2}-2}=\mathcal{Z}_{j},\qquad\psi_{j,\eta}:=\sqrt{\lambda
_{\eta}}M_{j}(x,\eta)
\]
then
\[
\Lambda_{j}f(x)=\sum_{\eta\in\mathcal{Z}_{j}} \langle
f,\psi_{j,\eta}\rangle\, \psi_{j,\eta}(x)\text{ }.
\]
\begin{proposition}
The family $(\psi_{j, \eta})_{j\in\mathbb{N},
\eta\in\mathcal{Z}_{j}}$ is a tight frame
\end{proposition}
\proof As
$$
\displaylines{
f = \lim_{J \to\infty}
\Big( L_{0}(f) + \sum_{j\leq J} \Lambda_{j}(f)\Big)\cr
\|f \|^{2} = \lim_{J \to\infty} \Big( \langle L_{0}(f),f \rangle
+\sum_{j\leq J} \langle\Lambda_{j}(f),f \rangle\Big)\cr
}
$$
but
\[
\langle\Lambda_{j}(f),f \rangle= \sum_{\eta\in\mathcal{Z}_{j}}
\langle f, \psi_{j, \eta} \rangle\langle\psi_{j, \eta},f \rangle=
\sum_{\eta\in
\mathcal{Z}_{j}} |\langle f, \psi_{j, \eta} \rangle|^{2}%
\]
and if $\psi_{0} $ is a normalized constant such that $\langle
L_{0}(f),f \rangle= |\langle f, \psi_{0} \rangle|^{2} $, then
\[
\|f\|^{2} = |\langle f, \psi_{0} \rangle|^{2} + \sum_{j \in\mathbb{N},
\eta\in\mathcal{Z}_{j}} |\langle f, \psi_{j, \eta} \rangle|^{2}\ .
\]
\hfill$\square$
\subsection{Localization properties}
Petrushev and coauthors have analysed the previous construction
proving that very nice localization properties hold.

In the case of the sphere of $\mathbb{R}^{d+1}$, where the spaces $H_{l}$ are
spanned by spherical harmonics, it is proved in
\cite{pnarco} the following localisation property:
for any $k$ there exists a constant $c_{k}$ such that :
\[
|\psi_{j\eta}(\xi)|\leq\frac{c_{k}2^{{d}\,\!{j}/2}}{(1+2^{j}\arccos
\langle\eta,\xi\rangle)^{k}%
}\text{ }\cdotp
\]
%
%
In the case of Jacobi polynomials on $[-1,1]$ with respect to the Jacobi
weight, it is
proved \cite{pxujacob} the following localization property:
for any $k$ there exist constants $C,c$ such that :
\[
|\psi_{j\eta}(\cos\theta)|\leq
\frac{c2^{j/2}}{(1+(2^{j}|\theta-\arccos
\eta|)^{k}\sqrt{w_{\alpha\beta}(2^{j},\cos\theta)}}\text{ }\cdotp
\]
where $w_{\alpha\beta}(n,x)=(1-x+n^{-2})^{\alpha+1/2}(1+x+n^{-2})^{\beta
+1/2},-1\leq x\leq1$ if $\alpha>-\frac 12,\ \beta>-\frac 12$. In the following
section, we consider the case, which is our framework,
where $\mathcal{Y}$ is the torus, $(e_{k})$ is the Fourier basis,
and $H_{m}=\hbox{Span}%
\{e_{m}\}$.
\subsection{ Quadrature Formula on the torus\label{}}
\begin{proposition} Assume that $\mathcal{Y}=\mathbb{T}$, the torus.  If, for $m\in\mathbb{N}$,
\[
\mathcal{K}_{m}=\Big\{\sum_{|k|\leq
m}a_{k}e^{ikx},a_{k}\in{\mathbb{C}}\Big\},
\]
the quadrature formula (\ref{quadr}), holds for
\[
\mathcal{X}_{m}=\Big\{
\frac{2l\pi}{m+1},\;l\in\{0,\ldots,m\}\Big\}
,\quad\lambda_{\frac{2l\pi}{m+1}}=\frac{1}{m+1}\cdotp
\]
\end{proposition}
\proof Let $\mathbb{T}$ be the torus, identified
with $[0,2\pi]$ and endowed with the measure
$\int_{\mathbb{T}}f\,d\mu=\frac{1}{2\pi}\int_{0}^{2\pi}f(x)\,dx$.
Let $f:\mathbb{T}\to\mathbb{C}$ with the following expansion on
the trigonometric basis,
\[
f=\sum_{k}a_{k}e^{ikx}\text{ }.
\]
For $m\in\mathbb{N}$ let us define,
\[
T_{m+1}(f)=\frac{1}{m+1}\sum_{l=0}^{m}f(x+\tfrac{2l\pi}{m+1})\text{
}.
\]
It is obvious that $T_{m+1}f$ is periodic with period
$\frac{2l\pi}{m+1}$. Therefore
$T_{m+1}f$ has the expansion
\begin{align*}
&\frac1{2\pi}\int_{\mathbb{T}}T_{m+1}f(x)\,e^{-ikx}\,dx=0~~\hbox{if}~~k\not=
0(\mathop{\rm mod}\;m+1)\text{ },\\
&\frac1{2\pi}\int_{\mathbb{T}}T_{m+1}(f)(x)\,e^{-ik(m+1)x}\,dx=a_{k{(m+1)}}\text{
}.
\end{align*}
Hence,
\[
T_{m+1}f(x)
=\sum
_{k}a_{k{(m+1)}}e^{ik(m+1)x}\text{ }.
\]
If $f$ is a trigonometric polynomial of degree smaller than or equal
to $m$, we have
\[
T_{m+1}f(x)=\frac{1}{m+1}\sum_{l=0}^{m}f(x+\tfrac{2l\pi}{m+1})=a_{0}%
=\frac1{2\pi}\int_{\mathbb{T}}f(u)du\text{ }.
\]
Therefore, if $f\in\mathcal{K}_{m}$,
\[
\frac{1}{m+1}\sum_{l=0}^{m}f(\tfrac{2l\pi
}{m+1})=\frac1{2\pi}\int_{\mathbb{T}}f(u)du\text{ .}%
\]
\qquad\hfill$\square$
\subsection{Localization properties for trigonometric series}
Following the steps of the sections above, we have for $\eta
\in\mathcal{X}_{2^{j+2}-1}$, that is
$\eta=\frac{2k\pi}{2^{j+2}},\;k\in\{0,\ldots ,2^{j+2}-1\}$,
\begin{equation}\label{psi}
\Lambda_j(x,y)=\sum_{l\neq0}a^{2}(\tfrac{l}{2^{j}})e^{il(x-y)}
,\qquad \psi_{j\eta}(x) =\frac1{2^{j/2+1}}
\sum_{l\neq0}a(\tfrac l{2^j})e^{il(x-\eta)}.
\end{equation}
%
\begin{figure}[h]
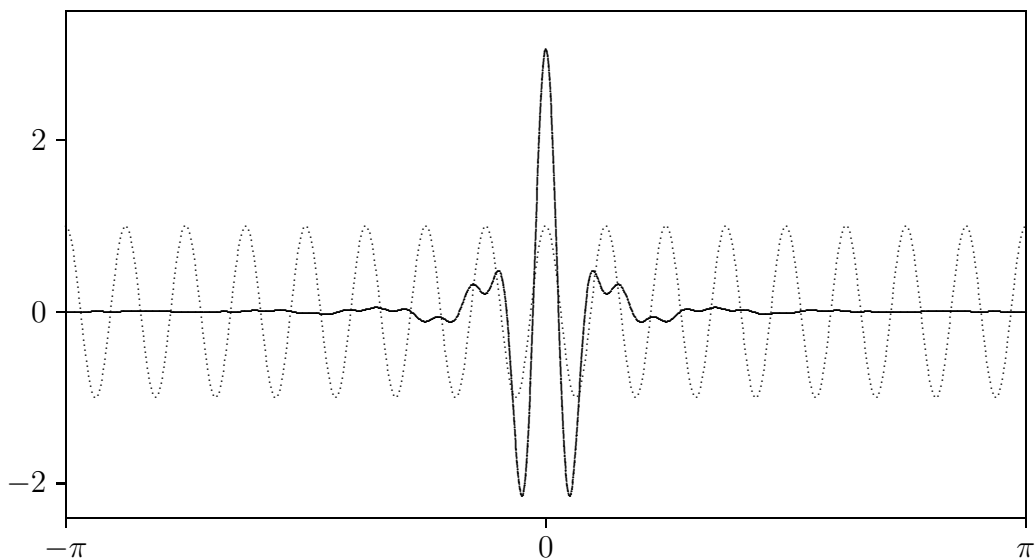

\hbox to\hsize\bgroup\hss
\beginpicture
\setcoordinatesystem units <.8truein,.45truein>
\setplotarea x from -3.1411 to 3.1411, y from -2.4 to 3.5
\axis bottom ticks short withvalues $-\pi$ $0$ $\pi$ / at -3.14 0 3.14 /  /
\axis top /
\axis right /
\axis left ticks short withvalues $-2$ $0$ $2$ / at -2 0 2 /  /
\setlinear
\plot
    -3.1415927    0.0000844
  -3.1353095    0.0000562
  -3.1290263  -0.0000273
  -3.1227431  -0.0001628
  -3.1164599  -0.0003450
  -3.1101767  -0.0005669
  -3.1038935  -0.0008195
  -3.0976104  -0.0010931
  -3.0913272  -0.0013765
  -3.085044   -0.0016586
  -3.0787608  -0.0019277
  -3.0724776  -0.0021730
  -3.0661944  -0.0023839
  -3.0599112  -0.0025514
  -3.0536281  -0.0026678
  -3.0473449  -0.0027273
  -3.0410617  -0.0027259
  -3.0347785  -0.0026622
  -3.0284953  -0.0025367
  -3.0222121  -0.0023522
  -3.0159289  -0.0021138
  -3.0096458  -0.0018286
  -3.0033626  -0.0015053
  -2.9970794  -0.0011542
  -2.9907962  -0.0007866
  -2.984513   -0.0004146
  -2.9782298  -0.0000504
  -2.9719467    0.0002938
  -2.9656635    0.0006065
  -2.9593803    0.0008772
  -2.9530971    0.0010970
  -2.9468139    0.0012586
  -2.9405307    0.0013566
  -2.9342475    0.0013881
  -2.9279644    0.0013525
  -2.9216812    0.0012514
  -2.915398     0.0010891
  -2.9091148    0.0008717
  -2.9028316    0.0006077
  -2.8965484    0.0003071
  -2.8902652  -0.0000186
  -2.8839821  -0.0003571
  -2.8776989  -0.0006954
  -2.8714157  -0.0010205
  -2.8651325  -0.0013197
  -2.8588493  -0.0015812
  -2.8525661  -0.0017946
  -2.8462829  -0.0019508
  -2.8399998  -0.0020430
  -2.8337166  -0.0020664
  -2.8274334  -0.0020188
  -2.8211502  -0.0019002
  -2.814867   -0.0017134
  -2.8085838  -0.0014633
  -2.8023006  -0.0011571
  -2.7960175  -0.0008041
  -2.7897343  -0.0004151
  -2.7834511  -0.0000019
  -2.7771679    0.0004226
  -2.7708847    0.0008454
  -2.7646015    0.0012537
  -2.7583183    0.0016351
  -2.7520352    0.0019787
  -2.745752     0.0022750
  -2.7394688    0.0025162
  -2.7331856    0.0026971
  -2.7269024    0.0028146
  -2.7206192    0.0028683
  -2.7143361    0.0028602
  -2.7080529    0.0027951
  -2.7017697    0.0026798
  -2.6954865    0.0025234
  -2.6892033    0.0023366
  -2.6829201    0.0021316
  -2.6766369    0.0019213
  -2.6703538    0.0017192
  -2.6640706    0.0015382
  -2.6577874    0.0013911
  -2.6515042    0.0012891
  -2.645221     0.0012419
  -2.6389378    0.0012570
  -2.6326546    0.0013398
  -2.6263715    0.0014926
  -2.6200883    0.0017153
  -2.6138051    0.0020045
  -2.6075219    0.0023544
  -2.6012387    0.0027561
  -2.5949555    0.0031986
  -2.5886723    0.0036687
  -2.5823892    0.0041514
  -2.576106     0.0046309
  -2.5698228    0.0050907
  -2.5635396    0.0055140
  -2.5572564    0.0058851
  -2.5509732    0.0061892
  -2.54469      0.0064132
  -2.5384069    0.0065463
  -2.5321237    0.0065803
  -2.5258405    0.0065099
  -2.5195573    0.0063332
  -2.5132741    0.0060513
  -2.5069909    0.0056689
  -2.5007078    0.0051935
  -2.4944246    0.0046358
  -2.4881414    0.0040089
  -2.4818582    0.0033282
  -2.475575     0.0026106
  -2.4692918    0.0018740
  -2.4630086    0.0011365
  -2.4567255    0.0004161
  -2.4504423  -0.0002700
  -2.4441591  -0.0009065
  -2.4378759  -0.0014799
  -2.4315927  -0.0019796
  -2.4253095  -0.0023976
  -2.4190263  -0.0027295
  -2.4127432  -0.0029743
  -2.40646    -0.0031344
  -2.4001768  -0.0032154
  -2.3938936  -0.0032264
  -2.3876104  -0.0031791
  -2.3813272  -0.0030875
  -2.375044   -0.0029674
  -2.3687609  -0.0028357
  -2.3624777  -0.0027099
  -2.3561945  -0.0026071
  -2.3499113  -0.0025435
  -2.3436281  -0.0025337
  -2.3373449  -0.0025899
  -2.3310617  -0.0027216
  -2.3247786  -0.0029350
  -2.3184954  -0.0032328
  -2.3122122  -0.0036137
  -2.305929   -0.0040727
  -2.2996458  -0.0046013
  -2.2933626  -0.0051872
  -2.2870795  -0.0058151
  -2.2807963  -0.0064671
  -2.2745131  -0.0071234
  -2.2682299  -0.0077630
  -2.2619467  -0.0083643
  -2.2556635  -0.0089060
  -2.2493803  -0.0093681
  -2.2430972  -0.0097323
  -2.236814   -0.0099833
  -2.2305308  -0.0101087
  -2.2242476  -0.0101001
  -2.2179644  -0.0099532
  -2.2116812  -0.0096685
  -2.205398   -0.0092505
  -2.1991149  -0.0087086
  -2.1928317  -0.0080562
  -2.1865485  -0.0073105
  -2.1802653  -0.0064918
  -2.1739821  -0.0056232
  -2.1676989  -0.0047293
  -2.1614157  -0.0038355
  -2.1551326  -0.0029674
  -2.1488494  -0.0021493
  -2.1425662  -0.0014037
  -2.136283   -0.0007503
  -2.1299998  -0.0002053
  -2.1237166    0.0002194
  -2.1174334    0.0005164
  -2.1111503    0.0006834
  -2.1048671    0.0007233
  -2.0985839    0.0006438
  -2.0923007    0.0004578
  -2.0860175    0.0001824
  -2.0797343  -0.0001614
  -2.0734512  -0.0005498
  -2.067168   -0.0009564
  -2.0608848  -0.0013539
  -2.0546016  -0.0017148
  -2.0483184  -0.0020124
  -2.0420352  -0.0022222
  -2.035752   -0.0023222
  -2.0294689  -0.0022946
  -2.0231857  -0.0021258
  -2.0169025  -0.0018077
  -2.0106193  -0.0013374
  -2.0043361  -0.0007179
  -1.9980529    0.0000420
  -1.9917697    0.0009282
  -1.9854866    0.0019215
  -1.9792034    0.0029984
  -1.9729202    0.0041317
  -1.966637     0.0052918
  -1.9603538    0.0064475
  -1.9540706    0.0075673
  -1.9477874    0.0086207
  -1.9415043    0.0095790
  -1.9352211    0.0104170
  -1.9289379    0.0111133
  -1.9226547    0.0116519
  -1.9163715    0.0120221
  -1.9100883    0.0122197
  -1.9038051    0.0122468
  -1.897522     0.0121119
  -1.8912388    0.0118297
  -1.8849556    0.0114206
  -1.8786724    0.0109096
  -1.8723892    0.0103261
  -1.866106     0.0097021
  -1.8598229    0.0090712
  -1.8535397    0.0084674
  -1.8472565    0.0079235
  -1.8409733    0.0074699
  -1.8346901    0.0071334
  -1.8284069    0.0069355
  -1.8221237    0.0068921
  -1.8158406    0.0070124
  -1.8095574    0.0072982
  -1.8032742    0.0077436
  -1.796991     0.0083354
  -1.7907078    0.0090530
  -1.7844246    0.0098691
  -1.7781414    0.0107506
  -1.7718583    0.0116596
  -1.7655751    0.0125547
  -1.7592919    0.0133927
  -1.7530087    0.0141297
  -1.7467255    0.0147232
  -1.7404423    0.0151335
  -1.7341591    0.0153253
  -1.727876     0.0152690
  -1.7215928    0.0149423
  -1.7153096    0.0143306
  -1.7090264    0.0134284
  -1.7027432    0.0122392
  -1.69646      0.0107756
  -1.6901768    0.0090596
  -1.6838937    0.0071214
  -1.6776105    0.0049985
  -1.6713273    0.0027350
  -1.6650441    0.0003796
  -1.6587609  -0.0020157
  -1.6524777  -0.0043978
  -1.6461946  -0.0067140
  -1.6399114  -0.0089140
  -1.6336282  -0.0109518
  -1.627345   -0.0127874
  -1.6210618  -0.0143882
  -1.6147786  -0.0157301
  -1.6084954  -0.0167987
  -1.6022123  -0.0175898
  -1.5959291  -0.0181092
  -1.5896459  -0.0183730
  -1.5833627  -0.0184065
  -1.5770795  -0.0182437
  -1.5707963  -0.0179256
  -1.5645131  -0.0174987
  -1.55823    -0.0170130
  -1.5519468  -0.0165202
  -1.5456636  -0.0160714
  -1.5393804  -0.0157149
  -1.5330972  -0.0154943
  -1.526814   -0.0154464
  -1.5205308  -0.0155994
  -1.5142477  -0.0159718
  -1.5079645  -0.0165711
  -1.5016813  -0.0173931
  -1.4953981  -0.0184222
  -1.4891149  -0.0196309
  -1.4828317  -0.0209812
  -1.4765485  -0.0224254
  -1.4702654  -0.0239080
  -1.4639822  -0.0253671
  -1.457699   -0.0267373
  -1.4514158  -0.0279516
  -1.4451326  -0.0289442
  -1.4388494  -0.0296532
  -1.4325663  -0.0300228
  -1.4262831  -0.0300061
  -1.4199999  -0.0295668
  -1.4137167  -0.0286811
  -1.4074335  -0.0273389
  -1.4011503  -0.0255448
  -1.3948671  -0.0233185
  -1.388584   -0.0206941
  -1.3823008  -0.0177197
  -1.3760176  -0.0144558
  -1.3697344  -0.0109737
  -1.3634512  -0.0073529
  -1.357168   -0.0036789
  -1.3508848  -0.0000395
  -1.3446017    0.0034774
  -1.3383185    0.0067877
  -1.3320353    0.0098137
  -1.3257521    0.0124875
  -1.3194689    0.0147531
  -1.3131857    0.0165690
  -1.3069025    0.0179100
  -1.3006194    0.0187678
  -1.2943362    0.0191524
  -1.288053     0.0190913
  -1.2817698    0.0186289
  -1.2754866    0.0178252
  -1.2692034    0.0167536
  -1.2629202    0.0154983
  -1.2566371    0.0141513
  -1.2503539    0.0128088
  -1.2440707    0.0115678
  -1.2377875    0.0105221
  -1.2315043    0.0097586
  -1.2252211    0.0093539
  -1.2189379    0.0093709
  -1.2126548    0.0098560
  -1.2063716    0.0108371
  -1.2000884    0.0123220
  -1.1938052    0.0142974
  -1.187522     0.0167290
  -1.1812388    0.0195623
  -1.1749557    0.0227242
  -1.1686725    0.0261255
  -1.1623893    0.0296633
  -1.1561061    0.0332254
  -1.1498229    0.0366938
  -1.1435397    0.0399493
  -1.1372565    0.0428761
  -1.1309734    0.0453660
  -1.1246902    0.0473232
  -1.118407     0.0486677
  -1.1121238    0.0493394
  -1.1058406    0.0493005
  -1.0995574    0.0485375
  -1.0932742    0.0470625
  -1.0869911    0.0449132
  -1.0807079    0.0421523
  -1.0744247    0.0388654
  -1.0681415    0.0351587
  -1.0618583    0.0311550
  -1.0555751    0.0269896
  -1.0492919    0.0228051
  -1.0430088    0.0187465
  -1.0367256    0.0149550
  -1.0304424    0.0115630
  -1.0241592    0.0086881
  -1.017876     0.0064283
  -1.0115928    0.0048575
  -1.0053096    0.0040217
  -0.9990265    0.0039366
  -0.9927433    0.0045855
  -0.9864601    0.0059191
  -0.9801769    0.0078564
  -0.9738937    0.0102862
  -0.9676105    0.0130709
  -0.9613274    0.0160502
  -0.9550442    0.0190471
  -0.9487610    0.0218732
  -0.9424778    0.0243362
  -0.9361946    0.0262465
  -0.9299114    0.0274247
  -0.9236282    0.0277086
  -0.9173451    0.0269602
  -0.9110619    0.0250711
  -0.9047787    0.0219683
  -0.8984955    0.0176178
  -0.8922123    0.0120274
  -0.8859291    0.0052477
  -0.8796459  -0.0026278
  -0.8733628  -0.0114650
  -0.8670796  -0.0210925
  -0.8607964  -0.0313065
  -0.8545132  -0.0418773
  -0.8482300  -0.0525564
  -0.8419468  -0.0630849
  -0.8356636  -0.0732026
  -0.8293805  -0.0826572
  -0.8230973  -0.0912133
  -0.8168141  -0.0986614
  -0.8105309  -0.1048261
  -0.8042477  -0.1095728
  -0.7979645  -0.1128140
  -0.7916813  -0.1145132
  -0.7853982  -0.1146870
  -0.7791150  -0.1134059
  -0.7728318  -0.1107928
  -0.7665486  -0.1070188
  -0.7602654  -0.1022981
  -0.7539822  -0.0968807
  -0.7476991  -0.0910429
  -0.7414159  -0.0850775
  -0.7351327  -0.0792823
  -0.7288495  -0.0739484
  -0.7225663  -0.0693480
  -0.7162831  -0.0657227
  -0.7099999  -0.0632726
  -0.7037168  -0.0621464
  -0.6974336  -0.0624331
  -0.6911504  -0.0641558
  -0.6848672  -0.0672674
  -0.6785840  -0.0716493
  -0.6723008  -0.0771123
  -0.6660176  -0.0834005
  -0.6597345  -0.0901975
  -0.6534513  -0.0971357
  -0.6471681  -0.1038073
  -0.6408849  -0.1097779
  -0.6346017  -0.1146012
  -0.6283185  -0.1178352
  -0.6220353  -0.1190582
  -0.6157522  -0.1178862
  -0.6094690  -0.1139880
  -0.6031858  -0.1071005
  -0.5969026  -0.0970414
  -0.5906194  -0.0837198
  -0.5843362  -0.0671440
  -0.5780530  -0.0474260
  -0.5717699  -0.0247827
  -0.5654867    0.0004671
  -0.5592035    0.0279085
  -0.5529203    0.0570409
  -0.5466371    0.0872917
  -0.5403539    0.1180325
  -0.5340708    0.1485988
  -0.5277876    0.1783112
  -0.5215044    0.2064984
  -0.5152212    0.2325204
  -0.5089380    0.2557919
  -0.5026548    0.2758050
  -0.4963716    0.2921492
  -0.4900885    0.3045295
  -0.4838053    0.3127808
  -0.4775221    0.3168784
  -0.4712389    0.3169435
  -0.4649557    0.3132444
  -0.4586725    0.3061913
  -0.4523893    0.2963266
  -0.4461062    0.2843095
  -0.4398230    0.2708955
  -0.4335398    0.2569115
  -0.4272566    0.2432267
  -0.4209734    0.2307207
  -0.4146902    0.2202488
  -0.4084070    0.2126067
  -0.4021239    0.2084950
  -0.3958407    0.2084848
  -0.3895575    0.2129862
  -0.3832743    0.2222196
  -0.3769911    0.2361930
  -0.3707079    0.2546839
  -0.3644247    0.2772289
  -0.3581416    0.3031201
  -0.3518584    0.3314106
  -0.3455752    0.3609267
  -0.3392920    0.3902901
  -0.3330088    0.4179469
  -0.3267256    0.4422057
  -0.3204425    0.4612810
  -0.3141593    0.4733442
  -0.3078761    0.4765782
  -0.3015929    0.4692357
  -0.2953097    0.4496993
  -0.2890265    0.4165416
  -0.2827433    0.3685838
  -0.2764602    0.3049496
  -0.2701770    0.2251158
  -0.2638938    0.1289535
  -0.2576106    0.0167627
  -0.2513274  -0.1107049
  -0.2450442  -0.2522325
  -0.2387610  -0.4061385
  -0.2324779  -0.5702887
  -0.2261947  -0.7421208
  -0.2199115  -0.9186832
  -0.2136283  -1.0966847
  -0.2073451  -1.2725561
  -0.2010619  -1.4425225
  -0.1947787  -1.6026838
  -0.1884956  -1.7491026
  -0.1822124  -1.8778965
  -0.1759292  -1.9853338
  -0.169646   -2.0679288
  -0.1633628  -2.1225355
  -0.1570796  -2.146436
  -0.1507964  -2.1374224
  -0.1445133  -2.0938688
  -0.1382301  -2.0147923
  -0.1319469  -1.8998999
  -0.1256637  -1.7496211
  -0.1193805  -1.5651239
  -0.1130973  -1.348314
  -0.1068142  -1.1018162
  -0.1005310  -0.8289393
  -0.0942478  -0.5336231
  -0.0879646  -0.2203701
  -0.0816814    0.1058374
  -0.0753982    0.4396330
  -0.0691150    0.7753730
  -0.0628319    1.1072512
  -0.0565487    1.4294188
  -0.0502655    1.7361082
  -0.0439823    2.021756
  -0.0376991    2.2811228
  -0.0314159    2.5094072
  -0.0251327    2.7023513
  -0.0188496    2.8563343
  -0.0125664    2.9684532
  -0.0062832    3.0365865
    0.           3.0594422
    0.0062832    3.0365865
    0.0125664    2.9684532
    0.0188496    2.8563343
    0.0251327    2.7023513
    0.0314159    2.5094072
    0.0376991    2.2811228
    0.0439823    2.021756
    0.0502655    1.7361082
    0.0565487    1.4294188
    0.0628319    1.1072512
    0.0691150    0.7753730
    0.0753982    0.4396330
    0.0816814    0.1058374
    0.0879646  -0.2203701
    0.0942478  -0.5336231
    0.1005310  -0.8289393
    0.1068142  -1.1018162
    0.1130973  -1.348314
    0.1193805  -1.5651239
    0.1256637  -1.7496211
    0.1319469  -1.8998999
    0.1382301  -2.0147923
    0.1445133  -2.0938688
    0.1507964  -2.1374224
    0.1570796  -2.146436
    0.1633628  -2.1225355
    0.169646   -2.0679288
    0.1759292  -1.9853338
    0.1822124  -1.8778965
    0.1884956  -1.7491026
    0.1947787  -1.6026838
    0.2010619  -1.4425225
    0.2073451  -1.2725561
    0.2136283  -1.0966847
    0.2199115  -0.9186832
    0.2261947  -0.7421208
    0.2324779  -0.5702887
    0.2387610  -0.4061385
    0.2450442  -0.2522325
    0.2513274  -0.1107049
    0.2576106    0.0167627
    0.2638938    0.1289535
    0.2701770    0.2251158
    0.2764602    0.3049496
    0.2827433    0.3685838
    0.2890265    0.4165416
    0.2953097    0.4496993
    0.3015929    0.4692357
    0.3078761    0.4765782
    0.3141593    0.4733442
    0.3204425    0.4612810
    0.3267256    0.4422057
    0.3330088    0.4179469
    0.3392920    0.3902901
    0.3455752    0.3609267
    0.3518584    0.3314106
    0.3581416    0.3031201
    0.3644247    0.2772289
    0.3707079    0.2546839
    0.3769911    0.2361930
    0.3832743    0.2222196
    0.3895575    0.2129862
    0.3958407    0.2084848
    0.4021239    0.2084950
    0.4084070    0.2126067
    0.4146902    0.2202488
    0.4209734    0.2307207
    0.4272566    0.2432267
    0.4335398    0.2569115
    0.4398230    0.2708955
    0.4461062    0.2843095
    0.4523893    0.2963266
    0.4586725    0.3061913
    0.4649557    0.3132444
    0.4712389    0.3169435
    0.4775221    0.3168784
    0.4838053    0.3127808
    0.4900885    0.3045295
    0.4963716    0.2921492
    0.5026548    0.2758050
    0.5089380    0.2557919
    0.5152212    0.2325204
    0.5215044    0.2064984
    0.5277876    0.1783112
    0.5340708    0.1485988
    0.5403539    0.1180325
    0.5466371    0.0872917
    0.5529203    0.0570409
    0.5592035    0.0279085
    0.5654867    0.0004671
    0.5717699  -0.0247827
    0.5780530  -0.0474260
    0.5843362  -0.0671440
    0.5906194  -0.0837198
    0.5969026  -0.0970414
    0.6031858  -0.1071005
    0.6094690  -0.1139880
    0.6157522  -0.1178862
    0.6220353  -0.1190582
    0.6283185  -0.1178352
    0.6346017  -0.1146012
    0.6408849  -0.1097779
    0.6471681  -0.1038073
    0.6534513  -0.0971357
    0.6597345  -0.0901975
    0.6660176  -0.0834005
    0.6723008  -0.0771123
    0.6785840  -0.0716493
    0.6848672  -0.0672674
    0.6911504  -0.0641558
    0.6974336  -0.0624331
    0.7037168  -0.0621464
    0.7099999  -0.0632726
    0.7162831  -0.0657227
    0.7225663  -0.0693480
    0.7288495  -0.0739484
    0.7351327  -0.0792823
    0.7414159  -0.0850775
    0.7476991  -0.0910429
    0.7539822  -0.0968807
    0.7602654  -0.1022981
    0.7665486  -0.1070188
    0.7728318  -0.1107928
    0.7791150  -0.1134059
    0.7853982  -0.1146870
    0.7916813  -0.1145132
    0.7979645  -0.1128140
    0.8042477  -0.1095728
    0.8105309  -0.1048261
    0.8168141  -0.0986614
    0.8230973  -0.0912133
    0.8293805  -0.0826572
    0.8356636  -0.0732026
    0.8419468  -0.0630849
    0.8482300  -0.0525564
    0.8545132  -0.0418773
    0.8607964  -0.0313065
    0.8670796  -0.0210925
    0.8733628  -0.0114650
    0.8796459  -0.0026278
    0.8859291    0.0052477
    0.8922123    0.0120274
    0.8984955    0.0176178
    0.9047787    0.0219683
    0.9110619    0.0250711
    0.9173451    0.0269602
    0.9236282    0.0277086
    0.9299114    0.0274247
    0.9361946    0.0262465
    0.9424778    0.0243362
    0.9487610    0.0218732
    0.9550442    0.0190471
    0.9613274    0.0160502
    0.9676105    0.0130709
    0.9738937    0.0102862
    0.9801769    0.0078564
    0.9864601    0.0059191
    0.9927433    0.0045855
    0.9990265    0.0039366
    1.0053096    0.0040217
    1.0115928    0.0048575
    1.017876     0.0064283
    1.0241592    0.0086881
    1.0304424    0.0115630
    1.0367256    0.0149550
    1.0430088    0.0187465
    1.0492919    0.0228051
    1.0555751    0.0269896
    1.0618583    0.0311550
    1.0681415    0.0351587
    1.0744247    0.0388654
    1.0807079    0.0421523
    1.0869911    0.0449132
    1.0932742    0.0470625
    1.0995574    0.0485375
    1.1058406    0.0493005
    1.1121238    0.0493394
    1.118407     0.0486677
    1.1246902    0.0473232
    1.1309734    0.0453660
    1.1372565    0.0428761
    1.1435397    0.0399493
    1.1498229    0.0366938
    1.1561061    0.0332254
    1.1623893    0.0296633
    1.1686725    0.0261255
    1.1749557    0.0227242
    1.1812388    0.0195623
    1.187522     0.0167290
    1.1938052    0.0142974
    1.2000884    0.0123220
    1.2063716    0.0108371
    1.2126548    0.0098560
    1.2189379    0.0093709
    1.2252211    0.0093539
    1.2315043    0.0097586
    1.2377875    0.0105221
    1.2440707    0.0115678
    1.2503539    0.0128088
    1.2566371    0.0141513
    1.2629202    0.0154983
    1.2692034    0.0167536
    1.2754866    0.0178252
    1.2817698    0.0186289
    1.288053     0.0190913
    1.2943362    0.0191524
    1.3006194    0.0187678
    1.3069025    0.0179100
    1.3131857    0.0165690
    1.3194689    0.0147531
    1.3257521    0.0124875
    1.3320353    0.0098137
    1.3383185    0.0067877
    1.3446017    0.0034774
    1.3508848  -0.0000395
    1.357168   -0.0036789
    1.3634512  -0.0073529
    1.3697344  -0.0109737
    1.3760176  -0.0144558
    1.3823008  -0.0177197
    1.388584   -0.0206941
    1.3948671  -0.0233185
    1.4011503  -0.0255448
    1.4074335  -0.0273389
    1.4137167  -0.0286811
    1.4199999  -0.0295668
    1.4262831  -0.0300061
    1.4325663  -0.0300228
    1.4388494  -0.0296532
    1.4451326  -0.0289442
    1.4514158  -0.0279516
    1.457699   -0.0267373
    1.4639822  -0.0253671
    1.4702654  -0.0239080
    1.4765485  -0.0224254
    1.4828317  -0.0209812
    1.4891149  -0.0196309
    1.4953981  -0.0184222
    1.5016813  -0.0173931
    1.5079645  -0.0165711
    1.5142477  -0.0159718
    1.5205308  -0.0155994
    1.526814   -0.0154464
    1.5330972  -0.0154943
    1.5393804  -0.0157149
    1.5456636  -0.0160714
    1.5519468  -0.0165202
    1.55823    -0.0170130
    1.5645131  -0.0174987
    1.5707963  -0.0179256
    1.5770795  -0.0182437
    1.5833627  -0.0184065
    1.5896459  -0.0183730
    1.5959291  -0.0181092
    1.6022123  -0.0175898
    1.6084954  -0.0167987
    1.6147786  -0.0157301
    1.6210618  -0.0143882
    1.627345   -0.0127874
    1.6336282  -0.0109518
    1.6399114  -0.0089140
    1.6461946  -0.0067140
    1.6524777  -0.0043978
    1.6587609  -0.0020157
    1.6650441    0.0003796
    1.6713273    0.0027350
    1.6776105    0.0049985
    1.6838937    0.0071214
    1.6901768    0.0090596
    1.69646      0.0107756
    1.7027432    0.0122392
    1.7090264    0.0134284
    1.7153096    0.0143306
    1.7215928    0.0149423
    1.727876     0.0152690
    1.7341591    0.0153253
    1.7404423    0.0151335
    1.7467255    0.0147232
    1.7530087    0.0141297
    1.7592919    0.0133927
    1.7655751    0.0125547
    1.7718583    0.0116596
    1.7781414    0.0107506
    1.7844246    0.0098691
    1.7907078    0.0090530
    1.796991     0.0083354
    1.8032742    0.0077436
    1.8095574    0.0072982
    1.8158406    0.0070124
    1.8221237    0.0068921
    1.8284069    0.0069355
    1.8346901    0.0071334
    1.8409733    0.0074699
    1.8472565    0.0079235
    1.8535397    0.0084674
    1.8598229    0.0090712
    1.866106     0.0097021
    1.8723892    0.0103261
    1.8786724    0.0109096
    1.8849556    0.0114206
    1.8912388    0.0118297
    1.897522     0.0121119
    1.9038051    0.0122468
    1.9100883    0.0122197
    1.9163715    0.0120221
    1.9226547    0.0116519
    1.9289379    0.0111133
    1.9352211    0.0104170
    1.9415043    0.0095790
    1.9477874    0.0086207
    1.9540706    0.0075673
    1.9603538    0.0064475
    1.966637     0.0052918
    1.9729202    0.0041317
    1.9792034    0.0029984
    1.9854866    0.0019215
    1.9917697    0.0009282
    1.9980529    0.0000420
    2.0043361  -0.0007179
    2.0106193  -0.0013374
    2.0169025  -0.0018077
    2.0231857  -0.0021258
    2.0294689  -0.0022946
    2.035752   -0.0023222
    2.0420352  -0.0022222
    2.0483184  -0.0020124
    2.0546016  -0.0017148
    2.0608848  -0.0013539
    2.067168   -0.0009564
    2.0734512  -0.0005498
    2.0797343  -0.0001614
    2.0860175    0.0001824
    2.0923007    0.0004578
    2.0985839    0.0006438
    2.1048671    0.0007233
    2.1111503    0.0006834
    2.1174334    0.0005164
    2.1237166    0.0002194
    2.1299998  -0.0002053
    2.136283   -0.0007503
    2.1425662  -0.0014037
    2.1488494  -0.0021493
    2.1551326  -0.0029674
    2.1614157  -0.0038355
    2.1676989  -0.0047293
    2.1739821  -0.0056232
    2.1802653  -0.0064918
    2.1865485  -0.0073105
    2.1928317  -0.0080562
    2.1991149  -0.0087086
    2.205398   -0.0092505
    2.2116812  -0.0096685
    2.2179644  -0.0099532
    2.2242476  -0.0101001
    2.2305308  -0.0101087
    2.236814   -0.0099833
    2.2430972  -0.0097323
    2.2493803  -0.0093681
    2.2556635  -0.0089060
    2.2619467  -0.0083643
    2.2682299  -0.0077630
    2.2745131  -0.0071234
    2.2807963  -0.0064671
    2.2870795  -0.0058151
    2.2933626  -0.0051872
    2.2996458  -0.0046013
    2.305929   -0.0040727
    2.3122122  -0.0036137
    2.3184954  -0.0032328
    2.3247786  -0.0029350
    2.3310617  -0.0027216
    2.3373449  -0.0025899
    2.3436281  -0.0025337
    2.3499113  -0.0025435
    2.3561945  -0.0026071
    2.3624777  -0.0027099
    2.3687609  -0.0028357
    2.375044   -0.0029674
    2.3813272  -0.0030875
    2.3876104  -0.0031791
    2.3938936  -0.0032264
    2.4001768  -0.0032154
    2.40646    -0.0031344
    2.4127432  -0.0029743
    2.4190263  -0.0027295
    2.4253095  -0.0023976
    2.4315927  -0.0019796
    2.4378759  -0.0014799
    2.4441591  -0.0009065
    2.4504423  -0.0002700
    2.4567255    0.0004161
    2.4630086    0.0011365
    2.4692918    0.0018740
    2.475575     0.0026106
    2.4818582    0.0033282
    2.4881414    0.0040089
    2.4944246    0.0046358
    2.5007078    0.0051935
    2.5069909    0.0056689
    2.5132741    0.0060513
    2.5195573    0.0063332
    2.5258405    0.0065099
    2.5321237    0.0065803
    2.5384069    0.0065463
    2.54469      0.0064132
    2.5509732    0.0061892
    2.5572564    0.0058851
    2.5635396    0.0055140
    2.5698228    0.0050907
    2.576106     0.0046309
    2.5823892    0.0041514
    2.5886723    0.0036687
    2.5949555    0.0031986
    2.6012387    0.0027561
    2.6075219    0.0023544
    2.6138051    0.0020045
    2.6200883    0.0017153
    2.6263715    0.0014926
    2.6326546    0.0013398
    2.6389378    0.0012570
    2.645221     0.0012419
    2.6515042    0.0012891
    2.6577874    0.0013911
    2.6640706    0.0015382
    2.6703538    0.0017192
    2.6766369    0.0019213
    2.6829201    0.0021316
    2.6892033    0.0023366
    2.6954865    0.0025234
    2.7017697    0.0026798
    2.7080529    0.0027951
    2.7143361    0.0028602
    2.7206192    0.0028683
    2.7269024    0.0028146
    2.7331856    0.0026971
    2.7394688    0.0025162
    2.745752     0.0022750
    2.7520352    0.0019787
    2.7583183    0.0016351
    2.7646015    0.0012537
    2.7708847    0.0008454
    2.7771679    0.0004226
    2.7834511  -0.0000019
    2.7897343  -0.0004151
    2.7960175  -0.0008041
    2.8023006  -0.0011571
    2.8085838  -0.0014633
    2.814867   -0.0017134
    2.8211502  -0.0019002
    2.8274334  -0.0020188
    2.8337166  -0.0020664
    2.8399998  -0.0020430
    2.8462829  -0.0019508
    2.8525661  -0.0017946
    2.8588493  -0.0015812
    2.8651325  -0.0013197
    2.8714157  -0.0010205
    2.8776989  -0.0006954
    2.8839821  -0.0003571
    2.8902652  -0.0000186
    2.8965484    0.0003071
    2.9028316    0.0006077
    2.9091148    0.0008717
    2.915398     0.0010891
    2.9216812    0.0012514
    2.9279644    0.0013525
    2.9342475    0.0013881
    2.9405307    0.0013566
    2.9468139    0.0012586
    2.9530971    0.0010970
    2.9593803    0.0008772
    2.9656635    0.0006065
    2.9719467    0.0002938
    2.9782298  -0.0000504
    2.984513   -0.0004146
    2.9907962  -0.0007866
    2.9970794  -0.0011542
    3.0033626  -0.0015053
    3.0096458  -0.0018286
    3.0159289  -0.0021138
    3.0222121  -0.0023522
    3.0284953  -0.0025367
    3.0347785  -0.0026622
    3.0410617  -0.0027259
    3.0473449  -0.0027273
    3.0536281  -0.0026678
    3.0599112  -0.0025514
    3.0661944  -0.0023839
    3.0724776  -0.0021730
    3.0787608  -0.0019277
    3.085044   -0.0016586
    3.0913272  -0.0013765
    3.0976104  -0.0010931
    3.1038935  -0.0008195
    3.1101767  -0.0005669
    3.1164599  -0.0003450
    3.1227431  -0.0001628
    3.1290263  -0.0000273
    3.1353095    0.0000562
    3.1415927    0.0000844
/
\setdots <2pt>
\plot
  -3.1415927    1.
  -3.1353095    0.9949510
  -3.1290263    0.9798551
  -3.1227431    0.9548645
  -3.1164599    0.9202318
  -3.1101767    0.8763067
  -3.1038935    0.8235326
  -3.0976104    0.7624425
  -3.0913272    0.6936533
  -3.085044     0.6178596
  -3.0787608    0.5358268
  -3.0724776    0.4483832
  -3.0661944    0.3564119
  -3.0599112    0.2608415
  -3.0536281    0.1626372
  -3.0473449    0.0627905
  -3.0410617  -0.0376902
  -3.0347785  -0.1377903
  -3.0284953  -0.236499
  -3.0222121  -0.3328195
  -3.0159289  -0.4257793
  -3.0096458  -0.5144395
  -3.0033626  -0.5979050
  -2.9970794  -0.6753328
  -2.9907962  -0.7459411
  -2.984513   -0.8090170
  -2.9782298  -0.8639234
  -2.9719467  -0.9101060
  -2.9656635  -0.9470983
  -2.9593803  -0.9745269
  -2.9530971  -0.9921147
  -2.9468139  -0.9996842
  -2.9405307  -0.9971589
  -2.9342475  -0.9845643
  -2.9279644  -0.9620277
  -2.9216812  -0.9297765
  -2.915398   -0.8881364
  -2.9091148  -0.8375280
  -2.9028316  -0.7784623
  -2.8965484  -0.7115357
  -2.8902652  -0.6374240
  -2.8839821  -0.5568756
  -2.8776989  -0.4707039
  -2.8714157  -0.3797791
  -2.8651325  -0.2850193
  -2.8588493  -0.1873813
  -2.8525661  -0.0878512
  -2.8462829    0.0125660
  -2.8399998    0.1128564
  -2.8337166    0.2120071
  -2.8274334    0.3090170
  -2.8211502    0.4029064
  -2.814867     0.4927273
  -2.8085838    0.5775727
  -2.8023006    0.6565858
  -2.7960175    0.7289686
  -2.7897343    0.7939904
  -2.7834511    0.8509945
  -2.7771679    0.8994053
  -2.7708847    0.9387339
  -2.7646015    0.9685832
  -2.7583183    0.9886517
  -2.7520352    0.9987370
  -2.745752     0.9987370
  -2.7394688    0.9886517
  -2.7331856    0.9685832
  -2.7269024    0.9387339
  -2.7206192    0.8994053
  -2.7143361    0.8509945
  -2.7080529    0.7939904
  -2.7017697    0.7289686
  -2.6954865    0.6565858
  -2.6892033    0.5775727
  -2.6829201    0.4927273
  -2.6766369    0.4029064
  -2.6703538    0.3090170
  -2.6640706    0.2120071
  -2.6577874    0.1128564
  -2.6515042    0.0125660
  -2.645221   -0.0878512
  -2.6389378  -0.1873813
  -2.6326546  -0.2850193
  -2.6263715  -0.3797791
  -2.6200883  -0.4707039
  -2.6138051  -0.5568756
  -2.6075219  -0.6374240
  -2.6012387  -0.7115357
  -2.5949555  -0.7784623
  -2.5886723  -0.8375280
  -2.5823892  -0.8881364
  -2.576106   -0.9297765
  -2.5698228  -0.9620277
  -2.5635396  -0.9845643
  -2.5572564  -0.9971589
  -2.5509732  -0.9996842
  -2.54469    -0.9921147
  -2.5384069  -0.9745269
  -2.5321237  -0.9470983
  -2.5258405  -0.9101060
  -2.5195573  -0.8639234
  -2.5132741  -0.8090170
  -2.5069909  -0.7459411
  -2.5007078  -0.6753328
  -2.4944246  -0.5979050
  -2.4881414  -0.5144395
  -2.4818582  -0.4257793
  -2.475575   -0.3328195
  -2.4692918  -0.236499
  -2.4630086  -0.1377903
  -2.4567255  -0.0376902
  -2.4504423    0.0627905
  -2.4441591    0.1626372
  -2.4378759    0.2608415
  -2.4315927    0.3564119
  -2.4253095    0.4483832
  -2.4190263    0.5358268
  -2.4127432    0.6178596
  -2.40646      0.6936533
  -2.4001768    0.7624425
  -2.3938936    0.8235326
  -2.3876104    0.8763067
  -2.3813272    0.9202318
  -2.375044     0.9548645
  -2.3687609    0.9798551
  -2.3624777    0.9949510
  -2.3561945    1.
  -2.3499113    0.9949510
  -2.3436281    0.9798551
  -2.3373449    0.9548645
  -2.3310617    0.9202318
  -2.3247786    0.8763067
  -2.3184954    0.8235326
  -2.3122122    0.7624425
  -2.305929     0.6936533
  -2.2996458    0.6178596
  -2.2933626    0.5358268
  -2.2870795    0.4483832
  -2.2807963    0.3564119
  -2.2745131    0.2608415
  -2.2682299    0.1626372
  -2.2619467    0.0627905
  -2.2556635  -0.0376902
  -2.2493803  -0.1377903
  -2.2430972  -0.236499
  -2.236814   -0.3328195
  -2.2305308  -0.4257793
  -2.2242476  -0.5144395
  -2.2179644  -0.5979050
  -2.2116812  -0.6753328
  -2.205398   -0.7459411
  -2.1991149  -0.8090170
  -2.1928317  -0.8639234
  -2.1865485  -0.9101060
  -2.1802653  -0.9470983
  -2.1739821  -0.9745269
  -2.1676989  -0.9921147
  -2.1614157  -0.9996842
  -2.1551326  -0.9971589
  -2.1488494  -0.9845643
  -2.1425662  -0.9620277
  -2.136283   -0.9297765
  -2.1299998  -0.8881364
  -2.1237166  -0.8375280
  -2.1174334  -0.7784623
  -2.1111503  -0.7115357
  -2.1048671  -0.6374240
  -2.0985839  -0.5568756
  -2.0923007  -0.4707039
  -2.0860175  -0.3797791
  -2.0797343  -0.2850193
  -2.0734512  -0.1873813
  -2.067168   -0.0878512
  -2.0608848    0.0125660
  -2.0546016    0.1128564
  -2.0483184    0.2120071
  -2.0420352    0.3090170
  -2.035752     0.4029064
  -2.0294689    0.4927273
  -2.0231857    0.5775727
  -2.0169025    0.6565858
  -2.0106193    0.7289686
  -2.0043361    0.7939904
  -1.9980529    0.8509945
  -1.9917697    0.8994053
  -1.9854866    0.9387339
  -1.9792034    0.9685832
  -1.9729202    0.9886517
  -1.966637     0.9987370
  -1.9603538    0.9987370
  -1.9540706    0.9886517
  -1.9477874    0.9685832
  -1.9415043    0.9387339
  -1.9352211    0.8994053
  -1.9289379    0.8509945
  -1.9226547    0.7939904
  -1.9163715    0.7289686
  -1.9100883    0.6565858
  -1.9038051    0.5775727
  -1.897522     0.4927273
  -1.8912388    0.4029064
  -1.8849556    0.3090170
  -1.8786724    0.2120071
  -1.8723892    0.1128564
  -1.866106     0.0125660
  -1.8598229  -0.0878512
  -1.8535397  -0.1873813
  -1.8472565  -0.2850193
  -1.8409733  -0.3797791
  -1.8346901  -0.4707039
  -1.8284069  -0.5568756
  -1.8221237  -0.6374240
  -1.8158406  -0.7115357
  -1.8095574  -0.7784623
  -1.8032742  -0.8375280
  -1.796991   -0.8881364
  -1.7907078  -0.9297765
  -1.7844246  -0.9620277
  -1.7781414  -0.9845643
  -1.7718583  -0.9971589
  -1.7655751  -0.9996842
  -1.7592919  -0.9921147
  -1.7530087  -0.9745269
  -1.7467255  -0.9470983
  -1.7404423  -0.9101060
  -1.7341591  -0.8639234
  -1.727876   -0.8090170
  -1.7215928  -0.7459411
  -1.7153096  -0.6753328
  -1.7090264  -0.5979050
  -1.7027432  -0.5144395
  -1.69646    -0.4257793
  -1.6901768  -0.3328195
  -1.6838937  -0.236499
  -1.6776105  -0.1377903
  -1.6713273  -0.0376902
  -1.6650441    0.0627905
  -1.6587609    0.1626372
  -1.6524777    0.2608415
  -1.6461946    0.3564119
  -1.6399114    0.4483832
  -1.6336282    0.5358268
  -1.627345     0.6178596
  -1.6210618    0.6936533
  -1.6147786    0.7624425
  -1.6084954    0.8235326
  -1.6022123    0.8763067
  -1.5959291    0.9202318
  -1.5896459    0.9548645
  -1.5833627    0.9798551
  -1.5770795    0.9949510
  -1.5707963    1.
  -1.5645131    0.9949510
  -1.55823      0.9798551
  -1.5519468    0.9548645
  -1.5456636    0.9202318
  -1.5393804    0.8763067
  -1.5330972    0.8235326
  -1.526814     0.7624425
  -1.5205308    0.6936533
  -1.5142477    0.6178596
  -1.5079645    0.5358268
  -1.5016813    0.4483832
  -1.4953981    0.3564119
  -1.4891149    0.2608415
  -1.4828317    0.1626372
  -1.4765485    0.0627905
  -1.4702654  -0.0376902
  -1.4639822  -0.1377903
  -1.457699   -0.236499
  -1.4514158  -0.3328195
  -1.4451326  -0.4257793
  -1.4388494  -0.5144395
  -1.4325663  -0.5979050
  -1.4262831  -0.6753328
  -1.4199999  -0.7459411
  -1.4137167  -0.8090170
  -1.4074335  -0.8639234
  -1.4011503  -0.9101060
  -1.3948671  -0.9470983
  -1.388584   -0.9745269
  -1.3823008  -0.9921147
  -1.3760176  -0.9996842
  -1.3697344  -0.9971589
  -1.3634512  -0.9845643
  -1.357168   -0.9620277
  -1.3508848  -0.9297765
  -1.3446017  -0.8881364
  -1.3383185  -0.8375280
  -1.3320353  -0.7784623
  -1.3257521  -0.7115357
  -1.3194689  -0.6374240
  -1.3131857  -0.5568756
  -1.3069025  -0.4707039
  -1.3006194  -0.3797791
  -1.2943362  -0.2850193
  -1.288053   -0.1873813
  -1.2817698  -0.0878512
  -1.2754866    0.0125660
  -1.2692034    0.1128564
  -1.2629202    0.2120071
  -1.2566371    0.3090170
  -1.2503539    0.4029064
  -1.2440707    0.4927273
  -1.2377875    0.5775727
  -1.2315043    0.6565858
  -1.2252211    0.7289686
  -1.2189379    0.7939904
  -1.2126548    0.8509945
  -1.2063716    0.8994053
  -1.2000884    0.9387339
  -1.1938052    0.9685832
  -1.187522     0.9886517
  -1.1812388    0.9987370
  -1.1749557    0.9987370
  -1.1686725    0.9886517
  -1.1623893    0.9685832
  -1.1561061    0.9387339
  -1.1498229    0.8994053
  -1.1435397    0.8509945
  -1.1372565    0.7939904
  -1.1309734    0.7289686
  -1.1246902    0.6565858
  -1.118407     0.5775727
  -1.1121238    0.4927273
  -1.1058406    0.4029064
  -1.0995574    0.3090170
  -1.0932742    0.2120071
  -1.0869911    0.1128564
  -1.0807079    0.0125660
  -1.0744247  -0.0878512
  -1.0681415  -0.1873813
  -1.0618583  -0.2850193
  -1.0555751  -0.3797791
  -1.0492919  -0.4707039
  -1.0430088  -0.5568756
  -1.0367256  -0.6374240
  -1.0304424  -0.7115357
  -1.0241592  -0.7784623
  -1.017876   -0.8375280
  -1.0115928  -0.8881364
  -1.0053096  -0.9297765
  -0.9990265  -0.9620277
  -0.9927433  -0.9845643
  -0.9864601  -0.9971589
  -0.9801769  -0.9996842
  -0.9738937  -0.9921147
  -0.9676105  -0.9745269
  -0.9613274  -0.9470983
  -0.9550442  -0.9101060
  -0.9487610  -0.8639234
  -0.9424778  -0.8090170
  -0.9361946  -0.7459411
  -0.9299114  -0.6753328
  -0.9236282  -0.5979050
  -0.9173451  -0.5144395
  -0.9110619  -0.4257793
  -0.9047787  -0.3328195
  -0.8984955  -0.236499
  -0.8922123  -0.1377903
  -0.8859291  -0.0376902
  -0.8796459    0.0627905
  -0.8733628    0.1626372
  -0.8670796    0.2608415
  -0.8607964    0.3564119
  -0.8545132    0.4483832
  -0.8482300    0.5358268
  -0.8419468    0.6178596
  -0.8356636    0.6936533
  -0.8293805    0.7624425
  -0.8230973    0.8235326
  -0.8168141    0.8763067
  -0.8105309    0.9202318
  -0.8042477    0.9548645
  -0.7979645    0.9798551
  -0.7916813    0.9949510
  -0.7853982    1.
  -0.7791150    0.9949510
  -0.7728318    0.9798551
  -0.7665486    0.9548645
  -0.7602654    0.9202318
  -0.7539822    0.8763067
  -0.7476991    0.8235326
  -0.7414159    0.7624425
  -0.7351327    0.6936533
  -0.7288495    0.6178596
  -0.7225663    0.5358268
  -0.7162831    0.4483832
  -0.7099999    0.3564119
  -0.7037168    0.2608415
  -0.6974336    0.1626372
  -0.6911504    0.0627905
  -0.6848672  -0.0376902
  -0.6785840  -0.1377903
  -0.6723008  -0.236499
  -0.6660176  -0.3328195
  -0.6597345  -0.4257793
  -0.6534513  -0.5144395
  -0.6471681  -0.5979050
  -0.6408849  -0.6753328
  -0.6346017  -0.7459411
  -0.6283185  -0.8090170
  -0.6220353  -0.8639234
  -0.6157522  -0.9101060
  -0.6094690  -0.9470983
  -0.6031858  -0.9745269
  -0.5969026  -0.9921147
  -0.5906194  -0.9996842
  -0.5843362  -0.9971589
  -0.5780530  -0.9845643
  -0.5717699  -0.9620277
  -0.5654867  -0.9297765
  -0.5592035  -0.8881364
  -0.5529203  -0.8375280
  -0.5466371  -0.7784623
  -0.5403539  -0.7115357
  -0.5340708  -0.6374240
  -0.5277876  -0.5568756
  -0.5215044  -0.4707039
  -0.5152212  -0.3797791
  -0.5089380  -0.2850193
  -0.5026548  -0.1873813
  -0.4963716  -0.0878512
  -0.4900885    0.0125660
  -0.4838053    0.1128564
  -0.4775221    0.2120071
  -0.4712389    0.3090170
  -0.4649557    0.4029064
  -0.4586725    0.4927273
  -0.4523893    0.5775727
  -0.4461062    0.6565858
  -0.4398230    0.7289686
  -0.4335398    0.7939904
  -0.4272566    0.8509945
  -0.4209734    0.8994053
  -0.4146902    0.9387339
  -0.4084070    0.9685832
  -0.4021239    0.9886517
  -0.3958407    0.9987370
  -0.3895575    0.9987370
  -0.3832743    0.9886517
  -0.3769911    0.9685832
  -0.3707079    0.9387339
  -0.3644247    0.8994053
  -0.3581416    0.8509945
  -0.3518584    0.7939904
  -0.3455752    0.7289686
  -0.3392920    0.6565858
  -0.3330088    0.5775727
  -0.3267256    0.4927273
  -0.3204425    0.4029064
  -0.3141593    0.3090170
  -0.3078761    0.2120071
  -0.3015929    0.1128564
  -0.2953097    0.0125660
  -0.2890265  -0.0878512
  -0.2827433  -0.1873813
  -0.2764602  -0.2850193
  -0.2701770  -0.3797791
  -0.2638938  -0.4707039
  -0.2576106  -0.5568756
  -0.2513274  -0.6374240
  -0.2450442  -0.7115357
  -0.2387610  -0.7784623
  -0.2324779  -0.8375280
  -0.2261947  -0.8881364
  -0.2199115  -0.9297765
  -0.2136283  -0.9620277
  -0.2073451  -0.9845643
  -0.2010619  -0.9971589
  -0.1947787  -0.9996842
  -0.1884956  -0.9921147
  -0.1822124  -0.9745269
  -0.1759292  -0.9470983
  -0.169646   -0.9101060
  -0.1633628  -0.8639234
  -0.1570796  -0.8090170
  -0.1507964  -0.7459411
  -0.1445133  -0.6753328
  -0.1382301  -0.5979050
  -0.1319469  -0.5144395
  -0.1256637  -0.4257793
  -0.1193805  -0.3328195
  -0.1130973  -0.236499
  -0.1068142  -0.1377903
  -0.1005310  -0.0376902
  -0.0942478    0.0627905
  -0.0879646    0.1626372
  -0.0816814    0.2608415
  -0.0753982    0.3564119
  -0.0691150    0.4483832
  -0.0628319    0.5358268
  -0.0565487    0.6178596
  -0.0502655    0.6936533
  -0.0439823    0.7624425
  -0.0376991    0.8235326
  -0.0314159    0.8763067
  -0.0251327    0.9202318
  -0.0188496    0.9548645
  -0.0125664    0.9798551
  -0.0062832    0.9949510
    0.           1.
    0.0062832    0.9949510
    0.0125664    0.9798551
    0.0188496    0.9548645
    0.0251327    0.9202318
    0.0314159    0.8763067
    0.0376991    0.8235326
    0.0439823    0.7624425
    0.0502655    0.6936533
    0.0565487    0.6178596
    0.0628319    0.5358268
    0.0691150    0.4483832
    0.0753982    0.3564119
    0.0816814    0.2608415
    0.0879646    0.1626372
    0.0942478    0.0627905
    0.1005310  -0.0376902
    0.1068142  -0.1377903
    0.1130973  -0.236499
    0.1193805  -0.3328195
    0.1256637  -0.4257793
    0.1319469  -0.5144395
    0.1382301  -0.5979050
    0.1445133  -0.6753328
    0.1507964  -0.7459411
    0.1570796  -0.8090170
    0.1633628  -0.8639234
    0.169646   -0.9101060
    0.1759292  -0.9470983
    0.1822124  -0.9745269
    0.1884956  -0.9921147
    0.1947787  -0.9996842
    0.2010619  -0.9971589
    0.2073451  -0.9845643
    0.2136283  -0.9620277
    0.2199115  -0.9297765
    0.2261947  -0.8881364
    0.2324779  -0.8375280
    0.2387610  -0.7784623
    0.2450442  -0.7115357
    0.2513274  -0.6374240
    0.2576106  -0.5568756
    0.2638938  -0.4707039
    0.2701770  -0.3797791
    0.2764602  -0.2850193
    0.2827433  -0.1873813
    0.2890265  -0.0878512
    0.2953097    0.0125660
    0.3015929    0.1128564
    0.3078761    0.2120071
    0.3141593    0.3090170
    0.3204425    0.4029064
    0.3267256    0.4927273
    0.3330088    0.5775727
    0.3392920    0.6565858
    0.3455752    0.7289686
    0.3518584    0.7939904
    0.3581416    0.8509945
    0.3644247    0.8994053
    0.3707079    0.9387339
    0.3769911    0.9685832
    0.3832743    0.9886517
    0.3895575    0.9987370
    0.3958407    0.9987370
    0.4021239    0.9886517
    0.4084070    0.9685832
    0.4146902    0.9387339
    0.4209734    0.8994053
    0.4272566    0.8509945
    0.4335398    0.7939904
    0.4398230    0.7289686
    0.4461062    0.6565858
    0.4523893    0.5775727
    0.4586725    0.4927273
    0.4649557    0.4029064
    0.4712389    0.3090170
    0.4775221    0.2120071
    0.4838053    0.1128564
    0.4900885    0.0125660
    0.4963716  -0.0878512
    0.5026548  -0.1873813
    0.5089380  -0.2850193
    0.5152212  -0.3797791
    0.5215044  -0.4707039
    0.5277876  -0.5568756
    0.5340708  -0.6374240
    0.5403539  -0.7115357
    0.5466371  -0.7784623
    0.5529203  -0.8375280
    0.5592035  -0.8881364
    0.5654867  -0.9297765
    0.5717699  -0.9620277
    0.5780530  -0.9845643
    0.5843362  -0.9971589
    0.5906194  -0.9996842
    0.5969026  -0.9921147
    0.6031858  -0.9745269
    0.6094690  -0.9470983
    0.6157522  -0.9101060
    0.6220353  -0.8639234
    0.6283185  -0.8090170
    0.6346017  -0.7459411
    0.6408849  -0.6753328
    0.6471681  -0.5979050
    0.6534513  -0.5144395
    0.6597345  -0.4257793
    0.6660176  -0.3328195
    0.6723008  -0.236499
    0.6785840  -0.1377903
    0.6848672  -0.0376902
    0.6911504    0.0627905
    0.6974336    0.1626372
    0.7037168    0.2608415
    0.7099999    0.3564119
    0.7162831    0.4483832
    0.7225663    0.5358268
    0.7288495    0.6178596
    0.7351327    0.6936533
    0.7414159    0.7624425
    0.7476991    0.8235326
    0.7539822    0.8763067
    0.7602654    0.9202318
    0.7665486    0.9548645
    0.7728318    0.9798551
    0.7791150    0.9949510
    0.7853982    1.
    0.7916813    0.9949510
    0.7979645    0.9798551
    0.8042477    0.9548645
    0.8105309    0.9202318
    0.8168141    0.8763067
    0.8230973    0.8235326
    0.8293805    0.7624425
    0.8356636    0.6936533
    0.8419468    0.6178596
    0.8482300    0.5358268
    0.8545132    0.4483832
    0.8607964    0.3564119
    0.8670796    0.2608415
    0.8733628    0.1626372
    0.8796459    0.0627905
    0.8859291  -0.0376902
    0.8922123  -0.1377903
    0.8984955  -0.236499
    0.9047787  -0.3328195
    0.9110619  -0.4257793
    0.9173451  -0.5144395
    0.9236282  -0.5979050
    0.9299114  -0.6753328
    0.9361946  -0.7459411
    0.9424778  -0.8090170
    0.9487610  -0.8639234
    0.9550442  -0.9101060
    0.9613274  -0.9470983
    0.9676105  -0.9745269
    0.9738937  -0.9921147
    0.9801769  -0.9996842
    0.9864601  -0.9971589
    0.9927433  -0.9845643
    0.9990265  -0.9620277
    1.0053096  -0.9297765
    1.0115928  -0.8881364
    1.017876   -0.8375280
    1.0241592  -0.7784623
    1.0304424  -0.7115357
    1.0367256  -0.6374240
    1.0430088  -0.5568756
    1.0492919  -0.4707039
    1.0555751  -0.3797791
    1.0618583  -0.2850193
    1.0681415  -0.1873813
    1.0744247  -0.0878512
    1.0807079    0.0125660
    1.0869911    0.1128564
    1.0932742    0.2120071
    1.0995574    0.3090170
    1.1058406    0.4029064
    1.1121238    0.4927273
    1.118407     0.5775727
    1.1246902    0.6565858
    1.1309734    0.7289686
    1.1372565    0.7939904
    1.1435397    0.8509945
    1.1498229    0.8994053
    1.1561061    0.9387339
    1.1623893    0.9685832
    1.1686725    0.9886517
    1.1749557    0.9987370
    1.1812388    0.9987370
    1.187522     0.9886517
    1.1938052    0.9685832
    1.2000884    0.9387339
    1.2063716    0.8994053
    1.2126548    0.8509945
    1.2189379    0.7939904
    1.2252211    0.7289686
    1.2315043    0.6565858
    1.2377875    0.5775727
    1.2440707    0.4927273
    1.2503539    0.4029064
    1.2566371    0.3090170
    1.2629202    0.2120071
    1.2692034    0.1128564
    1.2754866    0.0125660
    1.2817698  -0.0878512
    1.288053   -0.1873813
    1.2943362  -0.2850193
    1.3006194  -0.3797791
    1.3069025  -0.4707039
    1.3131857  -0.5568756
    1.3194689  -0.6374240
    1.3257521  -0.7115357
    1.3320353  -0.7784623
    1.3383185  -0.8375280
    1.3446017  -0.8881364
    1.3508848  -0.9297765
    1.357168   -0.9620277
    1.3634512  -0.9845643
    1.3697344  -0.9971589
    1.3760176  -0.9996842
    1.3823008  -0.9921147
    1.388584   -0.9745269
    1.3948671  -0.9470983
    1.4011503  -0.9101060
    1.4074335  -0.8639234
    1.4137167  -0.8090170
    1.4199999  -0.7459411
    1.4262831  -0.6753328
    1.4325663  -0.5979050
    1.4388494  -0.5144395
    1.4451326  -0.4257793
    1.4514158  -0.3328195
    1.457699   -0.236499
    1.4639822  -0.1377903
    1.4702654  -0.0376902
    1.4765485    0.0627905
    1.4828317    0.1626372
    1.4891149    0.2608415
    1.4953981    0.3564119
    1.5016813    0.4483832
    1.5079645    0.5358268
    1.5142477    0.6178596
    1.5205308    0.6936533
    1.526814     0.7624425
    1.5330972    0.8235326
    1.5393804    0.8763067
    1.5456636    0.9202318
    1.5519468    0.9548645
    1.55823      0.9798551
    1.5645131    0.9949510
    1.5707963    1.
    1.5770795    0.9949510
    1.5833627    0.9798551
    1.5896459    0.9548645
    1.5959291    0.9202318
    1.6022123    0.8763067
    1.6084954    0.8235326
    1.6147786    0.7624425
    1.6210618    0.6936533
    1.627345     0.6178596
    1.6336282    0.5358268
    1.6399114    0.4483832
    1.6461946    0.3564119
    1.6524777    0.2608415
    1.6587609    0.1626372
    1.6650441    0.0627905
    1.6713273  -0.0376902
    1.6776105  -0.1377903
    1.6838937  -0.236499
    1.6901768  -0.3328195
    1.69646    -0.4257793
    1.7027432  -0.5144395
    1.7090264  -0.5979050
    1.7153096  -0.6753328
    1.7215928  -0.7459411
    1.727876   -0.8090170
    1.7341591  -0.8639234
    1.7404423  -0.9101060
    1.7467255  -0.9470983
    1.7530087  -0.9745269
    1.7592919  -0.9921147
    1.7655751  -0.9996842
    1.7718583  -0.9971589
    1.7781414  -0.9845643
    1.7844246  -0.9620277
    1.7907078  -0.9297765
    1.796991   -0.8881364
    1.8032742  -0.8375280
    1.8095574  -0.7784623
    1.8158406  -0.7115357
    1.8221237  -0.6374240
    1.8284069  -0.5568756
    1.8346901  -0.4707039
    1.8409733  -0.3797791
    1.8472565  -0.2850193
    1.8535397  -0.1873813
    1.8598229  -0.0878512
    1.866106     0.0125660
    1.8723892    0.1128564
    1.8786724    0.2120071
    1.8849556    0.3090170
    1.8912388    0.4029064
    1.897522     0.4927273
    1.9038051    0.5775727
    1.9100883    0.6565858
    1.9163715    0.7289686
    1.9226547    0.7939904
    1.9289379    0.8509945
    1.9352211    0.8994053
    1.9415043    0.9387339
    1.9477874    0.9685832
    1.9540706    0.9886517
    1.9603538    0.9987370
    1.966637     0.9987370
    1.9729202    0.9886517
    1.9792034    0.9685832
    1.9854866    0.9387339
    1.9917697    0.8994053
    1.9980529    0.8509945
    2.0043361    0.7939904
    2.0106193    0.7289686
    2.0169025    0.6565858
    2.0231857    0.5775727
    2.0294689    0.4927273
    2.035752     0.4029064
    2.0420352    0.3090170
    2.0483184    0.2120071
    2.0546016    0.1128564
    2.0608848    0.0125660
    2.067168   -0.0878512
    2.0734512  -0.1873813
    2.0797343  -0.2850193
    2.0860175  -0.3797791
    2.0923007  -0.4707039
    2.0985839  -0.5568756
    2.1048671  -0.6374240
    2.1111503  -0.7115357
    2.1174334  -0.7784623
    2.1237166  -0.8375280
    2.1299998  -0.8881364
    2.136283   -0.9297765
    2.1425662  -0.9620277
    2.1488494  -0.9845643
    2.1551326  -0.9971589
    2.1614157  -0.9996842
    2.1676989  -0.9921147
    2.1739821  -0.9745269
    2.1802653  -0.9470983
    2.1865485  -0.9101060
    2.1928317  -0.8639234
    2.1991149  -0.8090170
    2.205398   -0.7459411
    2.2116812  -0.6753328
    2.2179644  -0.5979050
    2.2242476  -0.5144395
    2.2305308  -0.4257793
    2.236814   -0.3328195
    2.2430972  -0.236499
    2.2493803  -0.1377903
    2.2556635  -0.0376902
    2.2619467    0.0627905
    2.2682299    0.1626372
    2.2745131    0.2608415
    2.2807963    0.3564119
    2.2870795    0.4483832
    2.2933626    0.5358268
    2.2996458    0.6178596
    2.305929     0.6936533
    2.3122122    0.7624425
    2.3184954    0.8235326
    2.3247786    0.8763067
    2.3310617    0.9202318
    2.3373449    0.9548645
    2.3436281    0.9798551
    2.3499113    0.9949510
    2.3561945    1.
    2.3624777    0.9949510
    2.3687609    0.9798551
    2.375044     0.9548645
    2.3813272    0.9202318
    2.3876104    0.8763067
    2.3938936    0.8235326
    2.4001768    0.7624425
    2.40646      0.6936533
    2.4127432    0.6178596
    2.4190263    0.5358268
    2.4253095    0.4483832
    2.4315927    0.3564119
    2.4378759    0.2608415
    2.4441591    0.1626372
    2.4504423    0.0627905
    2.4567255  -0.0376902
    2.4630086  -0.1377903
    2.4692918  -0.236499
    2.475575   -0.3328195
    2.4818582  -0.4257793
    2.4881414  -0.5144395
    2.4944246  -0.5979050
    2.5007078  -0.6753328
    2.5069909  -0.7459411
    2.5132741  -0.8090170
    2.5195573  -0.8639234
    2.5258405  -0.9101060
    2.5321237  -0.9470983
    2.5384069  -0.9745269
    2.54469    -0.9921147
    2.5509732  -0.9996842
    2.5572564  -0.9971589
    2.5635396  -0.9845643
    2.5698228  -0.9620277
    2.576106   -0.9297765
    2.5823892  -0.8881364
    2.5886723  -0.8375280
    2.5949555  -0.7784623
    2.6012387  -0.7115357
    2.6075219  -0.6374240
    2.6138051  -0.5568756
    2.6200883  -0.4707039
    2.6263715  -0.3797791
    2.6326546  -0.2850193
    2.6389378  -0.1873813
    2.645221   -0.0878512
    2.6515042    0.0125660
    2.6577874    0.1128564
    2.6640706    0.2120071
    2.6703538    0.3090170
    2.6766369    0.4029064
    2.6829201    0.4927273
    2.6892033    0.5775727
    2.6954865    0.6565858
    2.7017697    0.7289686
    2.7080529    0.7939904
    2.7143361    0.8509945
    2.7206192    0.8994053
    2.7269024    0.9387339
    2.7331856    0.9685832
    2.7394688    0.9886517
    2.745752     0.9987370
    2.7520352    0.9987370
    2.7583183    0.9886517
    2.7646015    0.9685832
    2.7708847    0.9387339
    2.7771679    0.8994053
    2.7834511    0.8509945
    2.7897343    0.7939904
    2.7960175    0.7289686
    2.8023006    0.6565858
    2.8085838    0.5775727
    2.814867     0.4927273
    2.8211502    0.4029064
    2.8274334    0.3090170
    2.8337166    0.2120071
    2.8399998    0.1128564
    2.8462829    0.0125660
    2.8525661  -0.0878512
    2.8588493  -0.1873813
    2.8651325  -0.2850193
    2.8714157  -0.3797791
    2.8776989  -0.4707039
    2.8839821  -0.5568756
    2.8902652  -0.6374240
    2.8965484  -0.7115357
    2.9028316  -0.7784623
    2.9091148  -0.8375280
    2.915398   -0.8881364
    2.9216812  -0.9297765
    2.9279644  -0.9620277
    2.9342475  -0.9845643
    2.9405307  -0.9971589
    2.9468139  -0.9996842
    2.9530971  -0.9921147
    2.9593803  -0.9745269
    2.9656635  -0.9470983
    2.9719467  -0.9101060
    2.9782298  -0.8639234
    2.984513   -0.8090170
    2.9907962  -0.7459411
    2.9970794  -0.6753328
    3.0033626  -0.5979050
    3.0096458  -0.5144395
    3.0159289  -0.4257793
    3.0222121  -0.3328195
    3.0284953  -0.236499
    3.0347785  -0.1377903
    3.0410617  -0.0376902
    3.0473449    0.0627905
    3.0536281    0.1626372
    3.0599112    0.2608415
    3.0661944    0.3564119
    3.0724776    0.4483832
    3.0787608    0.5358268
    3.085044     0.6178596
    3.0913272    0.6936533
    3.0976104    0.7624425
    3.1038935    0.8235326
    3.1101767    0.8763067
    3.1164599    0.9202318
    3.1227431    0.9548645
    3.1290263    0.9798551
    3.1353095    0.9949510
    3.1415927    1.
/
\endpicture
\hss\egroup
\caption{Comparison between the wavelet $\psi_{j0}$ (solid) and
the trigonometric function at comparable frequency (dots). Here $j=4$.
The local character of the wavelet is apparent. In this picture the
wavelet has been normalized in order to have the same $L^2$ norm
(otherwise it would be smaller).\label{fig-wb}}
\end{figure}%
We prove  now at the same time the concentration of the wavelet
and of $\Lambda_{j}$. Let $A$ be a a continuous compactly
supported function and let us consider, associated to $A$ and $j$, the
function
\[
\xi_{j}(x)=\sum_{l\neq0}A(\tfrac{l}{2^{j}})e^{ilx}\ .
\]
We denote by $W_{1}^{k}=W_{1}^{k}(\mathbb{R})$ the
Sobolev space of functions with integrable weak derivative of order
$k.$
\begin{theorem}
\label{concen} Let $A$ be a continuous, compactly supported function
such that $A\in W_{1}^{k}(\mathbb{R})$ for some $k\geq2$.
Then, for all $m\in\mathbb{N}$ there exists a constant $c_{m,k}$ such that
\begin{equation}
|D^{m}\xi_{j}(x)|\leq\frac{c_{m,k}2^{(m+1)j}}{(1+|2^{j}x|)^{k}}\cdotp
 \label{conc}%
\end{equation}
\end{theorem}
\proof Clearly $\xi_{j}(x)$ is a trigonometric polynomial.
Moreover let us put
\[
B(x)=\bar{\mathcal{F}}(A)(x)=\frac{1}{2\pi}\int A(y)e^{iyx}\,dy.
\]
If $B\in\mathbb{L}_{1}(\mathbb{R}),$ by the Poisson summation formula :
\[
\xi_{j}(x)=\sum_{k}A(\tfrac{k}{2^{j}})e^{ikx}=\sum_{k}\hat{B}(\tfrac{k}{2^{j}%
})e^{ikx}=2\pi\sum_{l\in\mathbb{Z}}2^{j}B(2^{j}(x-2\pi l))
\]
and more generally if $D^{m}(B)\in\mathbb{L}_{1}(\mathbb{R})$
\[
D^{m}(\xi_{j})(x)=2\pi\sum_{l\in\mathbb{Z}}2^{(m+1)j}D^{m}(B)(2^{j}(x-2\pi
l))\text{ }.
\]
But
\[
D^{m}(B)(x)=D^{m}\bar{\mathcal{F}}(A)(x)=\bar{\mathcal{F}}((iy)^{m}A(y))(x)
\]
is bounded (and, by the Lebesgue-Riemann lemma, even vanishes at infinity). Furthermore,
$$
\displaylines{
(-ix)^{k}D^m(B)(x)=(-ix)^k\bar{\mathcal{F}}((iy)^mA(y))(x)=\bar
{\mathcal{F}}(D^{k}\{(iy)^{m}A(y)\})(x)=\cr
=i^{m}\sum_{l=0}^{k}\bar{\mathcal{F}}(\{D^{k-l}y^{m}D^{l}(A)(y)\})(x)\cr
}
$$
and this function is bounded as $A\in W_{1}^{k}$.
Therefore
\[
|D^{m}(B)(x)|\leq c_{m,k}\,\frac{1}{1+|x|^{k}}\leq c_{m,k}\,\frac{1}{(1+|x|)^{k}%
}\cdotp%
\]
Hence
\[
|D^{m}\xi_{j}(x)|\leq c_{m,k}2^{(m+1)j}\sum_{l\in\mathbb{Z}}\frac{1}%
{(1+|2^{j}(x-2\pi l)|)^{k}}\cdotp%
\]
The result is now a consequence of the following lemma.
\begin{lemma}
\label{ineq1}For $k\geq2$
\[
\theta_{j}(x) = \sum_{l\in\mathbb{Z}} \frac1{ (1+| 2^{j}
(x- 2\pi l)|)^{k} }%
\]
is a $2\pi$-periodic function such that
\[
\theta_{j}(x) \leq\frac5{ (1+| 2^{j} x|)^{k} }\cdotp%
\]
\end{lemma}
\proof Let $|t|\leq\pi$. Then
\[
\sum_{l\in\mathbb{Z}}\frac{1}{(1+|2^{j}(x-2\pi l)|)^{k}}=\frac{1}%
{(1+|2^{j}x|)^{k}}+\sum_{l\neq0}\frac{1}{(1+|2^{j}(x-2\pi l)|)^{k}}\cdotp%
\]%
Since
$$
\displaylines
{
\sum_{l\neq0}\frac{1}{(1+|2^{j}(x-2\pi l)|)^{k}}
 \leq2\sum_{l>0}\frac
{1}{(1+2^{j}(2l-1)\pi)^{k}}\le\cr
\leq\frac{2}{(1+2^{j}\pi)^{k}}\bigg(  1+\sum_{l\geq2}\bigg(  \frac
{1+2^{j}\pi}{1+2^{j}(2l-1)\pi}\bigg)  ^{k}\bigg)\le \cr
\leq\frac{2}{(1+2^{j}\pi)^{k}}\bigg(  1+\int_{2}^{\infty}\left(
\frac{1+2^{j}\pi}{1+2^{j}x\pi}\right)  ^{k}dx\bigg)
\leq\frac{4}{(1+2^{j}\pi)^{k}}\cr
}
$$
one gets finally
\[
\theta_{j}(x)\leq\frac{1}{(1+|2^{j}x|)^{k}}+\frac{4}{(1+2^{j}\pi)^{k}}%
\leq\frac{5}{(1+|2^{j}x|)^{k}}.
\]
\qquad\hfill$\square$
\begin{rem}
We are well aware that other choices of a well localized basis is
possible. As for instance the Lemari\'e-Meyer wavelet system
\cite{MR1228209}, that is also formed by trigonometric polynomials
and that is moreover orthonormal, whereas frames are not
orthonormal and redundant. Recall however that, as it is well
known, the advantage of redundant frames is robustness. Actually,
assume that $(e_i)_{i\in I}$ is an orthonormal basis of an Hilbert
space $H$, and that we have a noisy observation of a function
$f\in H$ by
\begin{equation}\label{err}
Y_i = \epsilon_i + \langle f, e_i\rangle
\end{equation}
where $(\epsilon_i)_i$ is a sequence of r.v.'s modeling the noise.
Then the error of the estimator $\sum_{i\in I} Y_i e_i$ is easily
computed
$$
\Big\| f- \sum_{i\in I} Y_i e_i \Big\|^2= \Big\| \sum_{i\in I}
\epsilon_i e_i \Big\|^2 = \sum_{i\in I} \epsilon_i^2\ .
$$
Conversely, if $(e_i)_{i\in I}$ was a tight frame then it is easy
to check the inequality
$$
  \sum_{i\in I}   |\langle \phi, e_i\rangle|^2 = \|\phi\|^2 =
  \Big\| \sum_{i\in I} \lambda_i e_i \Big\|^2
  \leq \sum_{i\in I} | \lambda_i|^2\ .
$$
Therefore the error is smaller
$$
\Big\| f- \sum_{i\in I} Y_i e_i \Big\|^2 = \Big\| \sum_{i\in I}
\epsilon_i e_i \Big\|^2 \leq \sum \epsilon_i^2\ .
$$
In fact often the inequality appears to be strict.
\end{rem}
\section{Assumptions and random wavelet coefficients}\label{assu}
\subsection{Assumptions on the model}
Consider the random field%
\[
X(\vartheta)=\sum_{l=-\infty}^{\infty}w_{l}e^{il\vartheta}\text{ ,
}\vartheta\in[0,2\pi],%
\]
where%
\[
w_{0}=0,\ Ew_{l}=0,\ E|w_{l}|^{2}=C_{l},\ %
l=1,2,\dots,\sum_{l=-\infty}^{\infty}C_{l}<\infty\text{ .}%
\]
Throughout this paper, we introduce the following regularity condition on the
behaviour of the angular power spectrum.
\medskip

\noindent\textbf{Assumption A1} There exists a function $g:\mathbb{N}%
\rightarrow\lbrack c_{1},c_{2}]$ such that $g\in W_{1}^{M}$ for some $M\geq0$
and%
\[
C_{l}=g(l)l^{-\alpha}\text{ for all }l\in\mathbb{N}\text{ , }\alpha>1\text{
}.
\]
For some results to follow, this assumption is strengthened to\medskip

\noindent
\textbf{Assumption A2} A1 holds and there exists a sequence of functions
$h_{N}(u):[\frac{1}{2},2]\rightarrow\lbrack\frac{c_{1}}{c_{2}},
\frac{c_{2}%
}{c_{1}}]$ such that
\[
h_{N}(\tfrac{4l}{N}):=\frac{g(l)}{g(\frac N4)},\qquad\frac{N}{8}\leq l\leq
\frac{N}{2}\text{ , }N=8,16,32,\dots\text{ }%
\]
and%
\begin{equation}
\sup_{N}\sup_{1/2\leq u\leq2}|h_{N}^{(M)}(u)|\leq C\text{ , some }C>0,\text{
some }M\in\mathbb{N}\text{ ,} \label{unbesov}%
\end{equation}
where $h_{N}^{(M)}$ denotes the $M$-th order weak derivative of
$h_N$.
\medskip

\noindent
\begin{rem}Condition A1 is a very mild regularity condition;
it is implied, for instance, if $g(l)$ is any trigonometric polynomial bounded
away from zero. The requirement $\alpha>1$ is necessary to ensure the sequence
$C_{l}$ to be summable, which in turn is a consequence of the finite variance of
the field. Condition A2 is a slightly stronger smoothness condition, which
implies that $h_{N}\in W_{1}^{M}$.
\end{rem}
\subsection{Random wavelet coefficients}

We recall the frame introduced in (\ref{psi}), namely
\[
\psi_{j\eta}(x)=\frac{1}{2^{1+j/2}}\sum_{l\neq0}a(\tfrac{l}{2^{j}
})e^{il(x-\eta)},\quad\eta=\frac{2k\pi}{2^{j+2}},\;k\in\{0,\ldots
,2^{j+2}-1\}.
\]
The notations will be shortened into the following way. For $j\in\mathbb{N}$,
we put $N=2^{j+2}$,
\begin{align*}
&\psi_{Nk}(t)   =\psi_{N}(t-k\tau),\quad\tau=\frac{2\pi}{N},\;k\in
\{0,\ldots,N-1\}\\
&\psi_{N}(t)    =\frac{1}{\sqrt{N}}\sum_{l=-\infty}^{\infty}{a}(\tfrac{4l}%
{N})e^{ilt}
\end{align*}
where $a$ is a $\mathcal{C}^{\infty}$ function, compactly supported in
$[\frac{1}{2},2]$. Hence, $\widehat{\psi}_{N}(l)$ has support in $(\frac{N}%
{8},\frac{N}{2})$; indeed%
\[
\widehat{\psi}_{N}(l)=\frac{1}{\sqrt{N}}\,a(\tfrac{4l}{N})\text{ .}%
\]
We have also
\[
\widehat{\psi}_{Nk}(l):=\frac 1{2\pi}\int_{-\pi}^{\pi}\psi_{N}%
(\vartheta+k\tau)e^{il\vartheta}\,d\vartheta=e^{-ilk\tau}
\widehat{\psi}_{N}(l)\text{ .}%
\]
We now define the associated wavelets coefficients of the process $X$%
\[
\beta_{Nk}:=\frac
1{2\pi}\int_{-\pi}^{\pi}X(\vartheta)\psi_{Nk}(\vartheta)\,d\vartheta
=\frac
1{2\pi}\int_{-\pi}^{\pi}X(\vartheta)\psi_{N}(\vartheta+k\tau)\,d\vartheta,
\quad k=0,1,\dots,N-1.%
\]
It is immediate to see that $E\beta_{Nk}=0;$ also
\begin{equation}\label{corr-coeff}
\begin{array}{c}
\displaystyle E(\beta_{Nk_{1}}\beta_{Nk_{2}})
=\int_{-\pi}^{\pi}\int_{-\pi}^{\pi
}EX(\vartheta)X(\vartheta^{\prime})\psi_{Nk_{1}}(\vartheta)\psi_{Nk_{2}%
}(\vartheta^{\prime})\,d\vartheta \,d\vartheta^{\prime}=\\
\displaystyle =\int_{-\pi}^{\pi}\int_{-\pi}^{\pi
}EX(\vartheta)X(\vartheta^{\prime})\psi_{Nk_{1}}(\vartheta)\psi_{Nk_{2}%
}(\vartheta^{\prime})\,d\vartheta \,d\vartheta^{\prime}=\cr
\displaystyle =\int_{-\pi}^{\pi}\int_{-\pi}^{\pi}\sum_{l=-\infty}^{\infty}C_{l}%
e^{il(\vartheta-\vartheta^{\prime})}\psi_{N}(\vartheta+k_{1}\tau)\psi
_{N}(\vartheta^{\prime}+k_{2}\tau)\,d\vartheta\,
d\vartheta^{\prime}=\\
\displaystyle =\sum_{l=-\infty}^{\infty}C_{l}\widehat{\psi}_{Nk_{1}}(l)\widehat{\psi
}_{Nk_{2}}(-l)=\sum_{l=-\infty}^{\infty}C_{l}|\widehat{\psi}_{N}(l)|^{2}%
e^{il\tau( k_{1}-k_{2})}=\\
\displaystyle=\frac 1N\sum_{l=-\infty}^{\infty}C_{l}a(\tfrac{4l}N)^2%
e^{il\tau( k_{1}-k_{2})}.
\end{array}
\end{equation}
%
Next we study the asymptotics  of the correlation coefficient, defined
as
\[
\Corr\left(  \beta_{Nk_{1}},\beta_{Nk_{2}}\right)
=\frac{\sum_{N/8\leq \left|  l\right|  \leq
N/2}C_{l}a^2(\tfrac{4l}{N})
e^{il\tfrac{2\pi}%
{N}(k_{1}-k_{2})}}{\sum_{N/8\leq|  l|  \leq N/2}C_{l}a^2%
(\frac{4l}{N})}\text{ .}%
\]
\begin{lemma}
\label{Lemma 1} Under Assumption A2,
\[
\left|  \Corr(  \beta_{Nk_{1}},\beta_{Nk_{2}})  \right|  \leq
\frac{2c_{M}}{(1+\left[  k_{1}-k_{2}\right]  _{N/2})^{M}}%
\]
for some $c_{M}>0$ where $[a]_{b}$ means $a(\operatorname{mod})b$.
\end{lemma}
\proof Write%
\[
\Corr(  \beta_{Nk_{1}},\beta_{Nk_{2}})  =\frac{\frac{1}{N}%
\sum_{N/8\leq|  l|  \leq N/2}\frac{C_{l}}{C_{N/4}}a^2(\frac
{l}{N})e^{il\tfrac{2\pi}{N}(k_{1}-k_{2})}}{\frac{1}{N}\sum_{N/8\leq\left|
l\right|  \leq N/2}\frac{C_{l}}{C_{N/4}}a^2(\frac{4l}{N})}\text{ ,}%
\]
and note that under Assumption A2 it holds for the denominator
\[
c_{1}2^{-\alpha}\leq\frac{1}{N}\sum_{N/8\leq\left|  l\right|  \leq N/2}%
\frac{C_{l}}{C_{N/4}}a^2(\tfrac{4l}{N})\leq c_{2}2^{\alpha}\text{
, as
}N\rightarrow\infty\text{ .}%
\]
Thus we focus on
$$
\displaylines{
\Big|  \frac{1}{N}\sum_{N/8\leq|  l|  \leq N/2}\frac{C_{l}%
}{C_{N/4}}a^2(\tfrac{4l}{N})e^{il\tfrac{2\pi}{N}(k_{1}-k_{2})}\Big|
\leq\cr
\le \frac{1}{N}c_{2}2^{\alpha}\Big|  \sum_{N/8\leq| l|
\leq
N/2}\frac{g(l)}{g(N/4)}(\tfrac{4l}{N})^{-\alpha}a^2(\tfrac{4l}{N}%
)\,e^{il\tfrac{2\pi}{N}(k_{1}-k_{2})}\Big|.\cr }
$$
To complete the argument, we use Theorem \ref{concen}. To apply the theorem,
we notice that for all $N$%
\[
A_{N}(\xi):=h_{N}(\xi)(\xi)^{-\alpha}a^2(\xi)
\]
is the (sampling of) the Fourier transform of an infinitely differentiable
expression, so we obtain%
$$
\displaylines{
\Big|  \frac 1N\sum_{N/8\leq|  l|  \leq N/2}\frac{g(l)}%
{g(N/4)}(\tfrac{4l}{N})^{-\alpha}a^2(\tfrac{4l}{N})e^{il\tfrac{2\pi}%
{N}(k_{1}-k_{2})}\Big| =\cr
=\Big| \frac 1N\sum_{N/8\leq|  l|  \leq N/2}h_{N}(\tfrac{4l}%
{N})(\tfrac{4l}{N})^{-\alpha}a^2(\tfrac{4l}{N})e^{il\tfrac{2\pi}{N}%
(k_{1}-k_{2})}\Big| \le\cr  \leq
\frac 1N\,\frac{2c_{M}N}{(1+N(\frac{2\pi}{N}\left[ k-k^{\prime}\right]
_{N/2}))^{M}}\leq\frac{2c_{M}}{(1+\left[  k-k^{\prime}\right]  _{N/2})^{M}%
}\text{ ,}%
} $$
which gives the required bound;
note that $c_{M}$ does not depend on $N,$ in view of
(\ref{unbesov}).

\hfill$\square$

\begin{rem} Lemma \ref{Lemma 1} highlights a quite remarkable property of random
wavelet coefficients. Indeed it entails that wavelet coefficients located at finite
distance are asymptotically (with respect to the frequency $N=2^{j+2}$)
uncorrelated.
\end{rem}
We write
\begin{align}
\sigma_N^2&:=\frac1N\sum_{{N/8\leq |  l|  \leq
N/2}}C_l a^2(\tfrac{4l}N)=\frac2N\sum_{{N/8\leq
l \leq N/2}}C_l a^2(\tfrac{4l}N)\label{varb}\\
\widehat{\beta}_{Nk}&:=\frac{\beta_{Nk}}{\sigma_{N}}
\raise2pt\hbox{.}\label{normb}
\end{align}
so that $E\widehat{\beta }_{Nk}^{2}=1$.
\smallskip

\begin{rem}\label{conv-var-coeff} It holds, as $N\to\infty$,
$$
\frac1N\sum_{{N/8\leq l \leq N/2}}a^2(\tfrac{4l}N)\enspace\to\enspace
\int_{1/2}^2a^2(t)\, dt
$$
therefore, as under Assumption A2 $0<c_1\le C_l/C_{N/4}\le
c_2<+\infty$,  there exist constants $0\le c'_1\le c'_2$ such that
\[
c'_1\leq\frac{\sigma_{N}^{2}}{C_{N/4}}\leq c'_2
\]
for every $N\ge 0$.
\end{rem}
In view of the asymptotic results of
next section, we shall always focus on the Gaussian case, as
motivated by our statistical applications.

\noindent\newline \textbf{Assumption B} The field
$X$ is Gaussian.

\section{Asymptotics of the wavelet statistics}\label{asym}

\subsection{The sample mean}

Our first aim in this Section is to investigate the asymptotic behaviour of
the sample mean for the random wavelet coefficients. More precisely, define%
\[
M_{N}:=\frac{1}{N}\sum_{k=1}^{N}\widehat{\beta}_{Nk}\text{ ;}%
\]
we have immediately $EM_{N}\equiv0$. Under Assumption A1, it is also simple to
show that
\begin{equation}\label{vanishing}
\begin
{array}{c}
\displaystyle E[M_N^2]=\frac1{N^2}\sum_{k_1,k_2=1}^{N}E[\widehat{\beta
}_{Nk_{1}}\widehat{\beta}_{Nk_{2}}]=\\
\displaystyle =\frac{1}{N^{2}}\sum_{N/8\leq|  l|  \leq N/2}\frac{C_{l}}%
{\sigma_{N}^{2}}a^2(\tfrac{4l}{N})\frac{1}{N}\Big|  \sum_{k=1}^{N}%
e^{i\tfrac{2\pi l}{N}k}\Big|  ^{2}=0,
\end{array}
\end{equation}
%
from the well-known properties of the Dirichlet kernel%
\[
\sum_{k=1}^{N}e^{i\tfrac{2\pi l}{N}k}=0\text{ for all }l\in\mathbb{N}\text{
such that }\tfrac{2\pi l}{N}\neq2k\pi\text{ , }k=0,\pm1,\pm2,\dots\text{ .}%
\]
It follows immediately that $M_N=0$  with
probability one. It is interesting to realize what happens here. Given the
fast decay of the covariances established in the previous Section, we might have
expected standard asymptotics to go through, i.e. a Central Limit Theorem for
the normalized sample mean. This is not the case because the variance is
degenerate; intuitively, this is due to the wavelet transform which is
'overdifferencing' the random field. Put in another way, if we
view the wavelet coefficients as a discrete time periodic random process, then
this process as a zero spectral density at the origin. This complicated
dependence structure does not prevent, however, the Central Limit Theorem to
hold for higher-order statistics, as we shall show in the sequel.
\subsection{Skewness and Kurtosis}\label{sk-ku-ssec}
Motivated by testing for non-Gaussianity on random fields on the torus, we
introduce here the Skewness and Kurtosis statistics of the wavelet coefficients.
More precisely, we consider (recall (\ref{varb}) and (\ref{normb}))%
\[
S_{N}:=\frac{1}{\sqrt{N}}\sum_{k=1}^{N}\widehat{\beta}_{Nk}^{3}
\quad\text{ and }\quad U_{N}:=\frac{1}{\sqrt{N}}\sum_{k=1}^{N}(  \widehat{\beta}_{Nk}%
^{4}-3)\ .%
\]
We have immediately%
\[
E[S_{N}]=E[U_{N}]=0\text{ }.
\]
The variance of these statistics is more complicated; in view of the following Lemma,
it is convenient to introduce Fej\'er's kernel, defined by%
$$
K_N(t):=\frac{1}{2\pi N}\Big|
\sum_{k=1}^{N}e^{itk}\Big|^2=\frac{1}{2\pi N}\,
\frac{\sin^2(\tfrac12\, Nt)}{\sin^2(\tfrac12\, t)}\cdotp
$$
Fej\'er's kernel vanishes at the
Fourier frequencies $\frac{2\pi}N\,l$, unless $l=kN$,
$k\in\mathbb{Z}$; moreover, if $l=kN,$ $K_{N}(\frac{2\pi}Nl)=
\frac N{2\pi}$. We define also%
\begin{equation}\label{varsk}
\sigma_{S_{N}}^{2}:=\frac{12\pi}{\sigma_{N}^{6}}\frac{1}{N^{3}}\sum
_{l_1l_2l_3}C_{l_1}a^2(\tfrac{4l_1}{N})C_{l_2}a^2%
(\tfrac{4l_2}{N})C_{l_3}a^2(\tfrac{4l_3}{N}) K_N(
(l_1+l_2+l_3)\tau), %
\end{equation}
and $\sigma_{U_{N}}^{2}=\sigma_{1U_{N}}^{2}+\sigma_{2U_{N}}^{2},$ where%
\begin{align}
&\sigma_{1U_{N}}^{2}:=\frac{72}{\sigma_{N}^{4}}\frac{2\pi}{N^{2}}\sum
_{l_1l_2}C_{l_1}a^2(\tfrac{4l_1}{N})C_{l_2}a^2(\tfrac{4l_2%
}{N})K_{N}((l_1+l_2)\tau)\text{ ,} \label{varkur1}\\
&\begin{matrix}\displaystyle
\sigma_{2U_N}^2:=\frac{24}{\sigma_{N}^{8}}\frac{48\pi}{N^4}
\sum_{l_1l_2l_3l_4}C_{l_1}a^2(\tfrac{4l_1}{N})C_{l_2}
a^2(\tfrac{4l_2}N)C_{l_3}a^2(\tfrac{4l_3}N)C_{l_4}a^2%
(\tfrac{4l_4}N)\times\\
\kern4cm\times K_N((l_1+l_2+l_3+l_4)\tau).
\end{matrix}\label{varkur2}%
\end{align}
\begin{rem} As Fej\'{e}r's kernel vanishes at the Fourier
frequencies $\frac{2\pi}N\l$, $l\not=kN$, $k\in\mathbb{Z}$, in the
previous expressions most terms vanish. Taking into account the
fact that $a^2(x)=0$ unless $\frac 12\le |x|\le 2$, one can write in
a more computationally tractable way
\begin{align*}
\sigma_{S_N}^{2} &  =\frac{6}{\sigma_{N}^{6}N^2}
\Big\{\sum_{l_1%
l_2}C_{l_1}a^2(\tfrac{4l_1}N)C_{l_2}a^2(\tfrac{4l_2}%
{N})C_{-l_1-l_2}a^2(-\tfrac{4l_1+4l_2}{N})+\\
&  +\sum_{l_1l_2}C_{l_1}{a}%
^{2}(\tfrac{4l_1}{N})C_{l_2}a^2(\tfrac{4l_2}{N})C_{N-l_1-l_2}%
a^2(4\tfrac{N-l_1-l_2}N)+\\
&  +\sum_{l_1l_2}C_{l_1}{a}%
^{2}(\tfrac{4l_1}N)C_{l_2}a^2(\tfrac{4l_2}{N})C_{-N-l_1-l_2}
a^2(-4\tfrac{N+l_1+l_2}{N})\Big\}%
\end{align*}
Likewise%
\[
\sigma_{1U_{N}}^{2}:=\frac{72}{\sigma_{N}^{4}}\frac{1}{N}\sum_{l_1}C_{l_1%
}^{2}{a}^{4}(\tfrac{4l_1}{N}).%
\]
A similar, a bit more complicated, expression can easily be derived
also for $\sigma_{2U_N}^2$.
\end{rem}
%
\begin{lemma}
Under Assumptions A1 and B%
\begin{align}
E[S_{N}^{2}]  &  =\sigma_{S_{N}}^{2}\text{ ,}\label{varskew}\\
E[U_{N}^{2}]  &  =\sigma_{1U_{N}}^{2}+\sigma_{2U_{N}}^{2}\text{ .}
\label{varkurt}%
\end{align}
\end{lemma}
\proof For (\ref{varskew}), we note that, by the diagram
formula (see \cite{MR611857}, p.108 e.g.)%
$$
\displaylines{
E[S_N^2]=\frac1{N\sigma_N^6}\sum_{k_1,k_2=1}^N%
E[\beta_{Nk_1}^3\beta_{Nk_2}^3]=\cr
=\frac1{N\sigma_N^6}\sum_{k_1,k_2=1}^N\big\{9
E[\beta_{Nk_1}\beta_{Nk_1}] E[\beta_{Nk_2}\beta_{Nk_2}]%
E[\beta_{Nk_1}\beta_{Nk_2}]  +6
E[\beta_{Nk_1}\beta_{Nk_2}]^3\big\}= \cr
=\frac{9}{N\sigma_{N}^{2}}\sum_{k_{1},k_{2}=1}^{N}E[\beta_{Nk_{1}}%
\beta_{Nk_{2}}]+\frac6{N\sigma_{N}^{2}}\sum_{k_1,k_2=1}^N
\Big\{\frac 1N
\sum_{l}C_{l}a^2(\tfrac{4l}{N})e^{il\tau( k_{1}-k_{2})}\Big\}
^{3}.\cr }
$$
We have seen already in (\ref{vanishing}) that the first term is
equal to zero. As for the second one, we obtain (recall that $\tau=\frac{2\pi}N$)%
$$
\displaylines{
\frac{6}{N\sigma_{N}^{6}}\sum_{k_{1},k_{2}=1}^{N}
\Big\{  \frac{1}{N}%
\sum_{l}C_{l}a^2(\tfrac{4l}{N})e^{il\tau( k_{1}-k_{2})}\Big\}
^{3}=\cr
=\frac{6}{N\sigma_{N}^{6}}\sum_{k_{1},k_{2}=1}^{N}\Big\{  \frac{1}{N^{3}%
}\sum_{l_1l_2l_3}C_{l_1}a^2(\tfrac{4l_1}{N})C_{l_2}{a}%
^{2}(\tfrac{4l_2}{N})C_{l_3}a^2(\tfrac{4l_3}{N}) e^{i\tau(
k_{1}-k_{2})(l_1+l_2+l_3)}\Big\} =\cr
=\frac{12\pi}{\sigma_{N}^{6}}\frac{1}{N^{3}}\sum_{l_1l_2l_3}C_{l_1%
}a^2(\tfrac{4l_1}{N})C_{l_2}a^2(\tfrac{4l_2}{N})C_{l_3}{a}%
^{2}(\tfrac{4l_3}{N})K_{N}((l_1+l_2+l_3)\tau)\text{ ,}%
}
$$
whence (\ref{varskew}) follows. Likewise, for (\ref{varkurt}), we
have
$$E[(\widehat \beta_{Nk_{1}}^{4}-3)(\widehat \beta_{Nk_2}^4-3)]=
E[\widehat \beta_{Nk_{1}}^{4}\widehat \beta_{Nk_2}^4]-9.
$$
Again by
the diagram formula and recalling that $E[\widehat \beta_{Nk_1}^2]=1$,
$$
E[\widehat \beta_{Nk_{1}}^{4}\widehat \beta_{Nk_2}^4]=
24E[\widehat \beta_{Nk_1}\widehat \beta_{Nk_2}]^4+
72E[\widehat \beta_{Nk_1}\widehat \beta_{Nk_2}]^2+9
$$
so that
$$
\displaylines{
EU_N^2  =\frac1N\sum_{k_1,k_2=1}^N
\big(E[\widehat \beta_{Nk_1}^4\widehat \beta_{Nk_2}^4]-9\big)=\cr
=\frac{24}{N\sigma_{N}^{8}}\sum_{k_1,k_2=1}^{N}E[\beta_{jk_{1}}%
\beta_{jk_2}]^4+\frac{72}{N\sigma_{N}^{4}}\sum_{k_1,k_2=1}^N%
E[\beta_{Nk_1}\beta_{Nk_2}]^2.\cr
}
$$
Now%
$$
\displaylines
{
\frac{1}{N\sigma_{N}^{8}}\sum_{k_{1},k_{2}=1}^{N}
E[\beta_{Nk_1}\beta_{Nk_2}]^{4}=\frac1{N^5\sigma_N^8}
\sum_{k_1,k_2=1}^N\Big(\sum_{l}C_{l}a^2(\tfrac{4l}{N})e^{il\tau(k_{1}-k_{2})}
\Big)^{4}=\cr
=\frac{1}{N^5\sigma_{N}^{8}}\!\sum_{k_{1},k_{2}=1}^{N}%
\sum_{l_1l_2l_3l_{4}}\!\!\!C_{l_1}a^2(\tfrac{4l_1}{N})C_{l_2}{a}%
^{2}(\tfrac{4l_2}{N})C_{l_3}a^2(\tfrac{4l_3}{N})C_{l_{4}}a^2%
(\tfrac{4l_{4}}{N})e^{i\tau(k_1-k_2)(l_1+l_2+l_3+l_{4})}=\cr
=\frac{2\pi}{N^4\sigma_N^8}\sum_{l_1l_2l_3l_{4}%
}C_{l_1}a^2(\tfrac{4l_1}{N})C_{l_2}a^2(\tfrac{4l_2}%
{N})C_{l_3}a^2(\tfrac{4l_3}{N})C_{l_{4}}a^2(\tfrac{4l_{4}}{N}%
)K_{N}((l_1+l_2+l_3+l_{4})\tau)\ .\cr }
$$
On the other hand
$$
\displaylines
{
\frac1{N\sigma_{N}^4}\sum_{k_1,k_2=1}^NE[\beta_{Nk_{1}}\beta
_{Nk_2}]^2  =\frac1{N^3\sigma_{N}^{4}}\sum_{k_{1}%
,k_{2}=1}^{N}\Big(\sum_{l}C_{l}a^2(\tfrac{4l}{N})
e^{i\tau l(k_{1}-k_{2})}\Big)^{2}=\cr
=\frac{2\pi}{N^2\sigma_N^4}\sum_{l_1l_2}C_{l_1}%
a^2(\tfrac{4l_1}{N})C_{l_2}a^2(\tfrac{4l_2}{N})K_{N}((l_1+l_2)\tau)\ .\cr }
$$
\hfill$\square$
\begin{rem}\label{conv-var-skku} By the same arguments of Remark
\ref{conv-var-coeff} it is immediate that, as $N\to\infty$,
$$
\displaylines{
\frac1{N^3}\sum_{l_1l_2l_3}a^2(\tfrac{4l_1}N)a^2%
(\tfrac{4l_2}N)a^2(\tfrac{4l_3}N) K_N(\tfrac{2\pi}N
(l_1+l_2+l_3))\enspace \to\enspace\cr
\int_{1/2}^2\int_{1/2}^2\int_{1/2}^2
a^2(t_1)a^2(t_2)a^2(t_3)\,{t_1+t_2+t_3\over
2\sin^2(t_1+t_2+t_3)}\,dt_1\,dt_2\,dt_3\ .\cr
}
$$
Therefore under Assumptions A1 and B we have
\begin{equation}\label{bound-var-skku}
0<c_1\le\sigma_{S_N}^2\le c_2
\end{equation}
for some constants $0\le c_1\le c_2$. By the same argument we see that also
$\sigma_{U_N}^2$ is bounded and bounded away from $0$, for every
$N>0$
\end{rem}
\section{The Central Limit Theorem}\label{tcl}
This Section is devoted to the Central Limit Theorem for our statistics of
interest. The idea is to prove the results by  the method of moments.
To analyze the behaviour of higher order moments, we shall extensively
use the already mentioned diagram formula, that states that, for a
multivariate zero-mean Gaussian vector $(X_{1},\dots,X_{2k}),$ it holds
\begin{equation}
E(X_1 X_2 \ldots  X_{2k})=\sum E(X_{i_1}X_{i_2}%
)\dots E(X_{i_{2k-1}}X_{i_{2k}})\ . \label{adle}%
\end{equation}
Consider the cartesian product $I\times J,$ where $I,J$ are sets of
positive integers of cardinality $\#(I\mathcal{)}=P,$ $\#(J\mathcal{)}=Q;$
it is convenient to visualize these elements in a $P\times Q$ matrix
with $P$
rows. A {\it diagram} $\gamma
$ is any partition of the $P\cdot Q$ elements into pairs like
$\{(i_1,j_1),(i_2,j_2)\}$; these pairs are called the {\it edges}
of the diagram. We label $\Gamma (I,J)$ the family of these diagrams.
It can be checked that, for given $I,J,$ there exist
$(P\cdot Q-1)!!$
different diagrams, each of them composed by $\frac 12\,P\cdot Q$ pairs; we
recall that $(2p-1)!!:=(2p-1)\cdot (2p-3)\cdots 1$
for $p=1,2,\dots$

To any diagram we can associate a
{\it graph} having $I$ as the set of vertices (or nodes) and
connecting any two of these vertices, $i_k, i_{k'}$, by an arc
whenever in the diagram an edge of the type $\{ (i_k,j_k),
(i_{k'},j_{k'})\}$ is present. This graph is not directed, that is,
$(i_1,i_2)$ and $(i_2,i_1)$ identify the
same arc; however, we do allow for repetitions of edges
(two rows may be
linked twice).
We shall use some
result on graphs below; with a slight abuse of notation, we denote the graph
$\gamma $ with the same letter as the corresponding diagram.

We say that

a) A diagram has a \emph{flat edge} if there is at least a pair
$\{(i_1,j_1),(i_2,j_2)\}$ with $i_1=i_2$. we write $%
\gamma\in\Gamma_{F}(I,J)$ for a diagram with at least a flat edge, and $%
\gamma\in\Gamma_{\overline{F}}(I,J)$ otherwise. A graph corresponding to a
diagram with a flat edge includes an edge $i_{k}i_{k}$ which arrives in the
same vertex where it started.

b) A diagram $\gamma\in\Gamma_{\overline{F}}(I,J)$ is \emph{connected} if it
is not possible to partition the corresponding graph
into two sets $A,B$ such that there
are no arcs connecting the nodes in $A$ with nodes in $B$. We write $\gamma\in\Gamma
_{C}(I,J)$ for connected diagrams, $\gamma\in\Gamma_{\overline{C}}(I,J)$
otherwise.

c) A diagram $\gamma\in\Gamma_{\overline{F}}(I,J)$ is \emph{paired}
if,
given any two sets of edges $\{(i_1,j_1),(i_2,j_2)\} $ and $\{
(i_3,j_3),(i_4,j_4)\}$, then $i_1=i_3$ implies
$i_2=i_4$; in words, the rows are completely coupled two by two. We write
$\gamma\in\Gamma_{P}(I,J)$ for paired diagrams.

Obviously if $P>2$ a paired diagram cannot be connected. Note that
if $Q$ is odd, paired diagrams cannot have flat edges, so that the
assumption $\gamma\in\Gamma_{\overline{F}}(I,J)$ becomes redundant.
If $I$ or $J$ (or both) can be simply taken as
the set of the first $p$ or $q$ natural numbers, i.e. $I=\{
1,\dots,p\}$, $J=\{  1,\dots,q\}$
we shall occasionally write $\Gamma(I,q),\Gamma
p,J)$ or $\Gamma(p,q)$ for $\Gamma(I,J).$
For a nice and comprehensive discussion on the diagram formula, we refer to
\cite{MR1956046}; see also \cite{MR2118863} for recent advances in this area.
\begin{theorem}
\label{TLC}Under Assumptions A1 and B, as $N\rightarrow\infty,$
\[
\lim_{N\rightarrow\infty}E\Big[\Big(  \frac{S_{N}}{\sigma_{S_{N}}}
\Big)
^{p_{1}}\Big(   \frac{U_{N}}{\sigma_{U_{N}}}\Big)  ^{p_{2}}\Big]=
\begin{cases}
(2p_1-1)!!\cdot (2p_2-1)!!&\text{ if }p_{1},p_{2}
\text{ are even}\\
0&\text{ otherwise.}%
\end{cases}
%
\]
Hence%
\[
\begin{pmatrix}
\vrule height0pt depth12pt width0pt \frac 1{\sigma_{S_N}}S_N\\
\frac 1{\sigma_{U_N}}U_N%
\end{pmatrix}
\enspace\mathop{\rightarrow}^{\mathcal D}_{N\to\infty}\enspace
N(0,I_2)\ .
%
\]
\end{theorem}
\proof For brevity's sake and notational simplicity, we focus
on the case $p_2=0;$ the general argument can be pursued under the
same lines. The idea is to use (\ref{adle}) and partition the various summands
in this expression according to the nature of the associated diagrams/graphs.
To this aim, we visualize our coefficients as positioned on a diagram with
$I=2p_{1}=2p$ rows and $J=3$ columns; we obtain%
$$
\displaylines{ ES_{N}^{2p}=\frac{1}{N^{p}\sigma_{N}^{6p}}
\sum_{k_1,\dots, k_{2p}}E
[
\beta_{Nk_1}^3\dots \beta_{Nk_{2p}}^{3}
]=\cr
=\frac{1}{N^{p}\sigma_{N}^{6p}}\sum_{k_1,\dots, k_{2p}}E\Big[
\prod_{k\in\{k_1,\dots, k_{2p}\}}\Big\{  \frac{1}{\sqrt{N}}\sum_{N/8\leq|
l|  \leq N/2}w_l\,a(\tfrac{4l}N)\,e^{il\tau
k}\Big\}^3\Big] =\cr
=\frac1{N^{4p}\sigma_N^{6p}}\sum_{k_1, \dots, k_{2p}}E\Big[
\prod_{k\in\{k_1,\dots ,k_{2p}\}}\prod_{j=1}^{3}\sum_{N/8\leq|  l|
\leq N/2}w_l\,a(\tfrac{4l}N)e^{il\tau
k}\Big]=\cr
=\frac{1}{N^{4p}\sigma_{N}^{6p}}\sum_{k_1,\dots,
k_{2p}}\sum_{\gamma\in \Gamma(2p,3)}\prod_{\{  (m,j)(m^{\prime},j^{\prime})\}  \in\gamma}
\!\Big\{\!\sum_{N/8\le|l|,|l|'\le N/2}\!\!\!\!\!E[w_lw_{l'}]\,
a(\tfrac{4l}N)
\, a(\tfrac{4l'}N)e^{il\tau k_m}e^{il'\tau
k_{m'}}\Big\} \cr
=\frac{1}{N^{4p}\sigma_{N}^{6p}}\sum_{k_1,\dots ,k_{2p}}\sum_{\gamma\in
\Gamma(2p,3)}\prod_{\{  (m,j)(m^{\prime},j^{\prime})\}  \in\gamma
}\Big\{  \sum_{N/8\le |l|\le N/2}C_l\, a^2(\tfrac{4l}N)\,e^{il\tau k_m}
e^{-il\tau k_{m'}}\Big\}\ ,\cr%
}
$$
where we used the fact that $E[w_lw_{l'}]=0$ unless $l'=-l$. For
fixed $k_1,\dots, k_{2p}$, let us concentrate on the contribution
of a single diagram $\gamma\in\Gamma(2p,3)$.
$$
\displaylines{ \prod_{\{  (m,j)(m^{\prime},j^{\prime})\}  \in\gamma
}\Big\{  \sum_{N/8\le |l|\le N/2}C_l\, a^2(\tfrac{4l}N)\,e^{il\tau
k_m} e^{-il\tau k_{m'}}\Big\}=\cr \sum_{l_{11},\dots,
l_{13},\dots,l_{2p,1},\dots,l_{2p,3}} \prod_{\{
(m,j)(m^{\prime},j^{\prime})\}  \in\gamma
}C_{l_{m,j}}a^2(\tfrac{4l_{m,j}}N)e^{il_{m,j}\tau k_m} e^{-il_{m,j}\tau
k_{m'}}\,\delta_{l_{m,j},l_{m',j'}}\cr }
$$
where the indices $l_{m,j}$, $m=1,\dots, 2p$, $j=1,2,3$ vary between
$N/8$ and $N/2$. Now, as every $k_m$ appears exactly three times in the diagram,
$$
\displaylines{ \prod_{\{  (m,j)(m^{\prime},j^{\prime})\}  \in\gamma
}C_{l_{m,j}}a^2(\tfrac{4l_{m,j}}N)e^{il_{m,j}\tau k_m} e^{-il_{m,j}\tau
k_{m'}}\,\delta_{l_{m,j},l_{m',j'}}=\cr
=\Big(\prod_{\{
(m,j)(m^{\prime},j^{\prime})\}\in\gamma
}C_{l_{m,j}}a^2(\tfrac{4l_{m,j}}N)\,\delta_{l_{m,j},l_{m',j'}}\Big)\Big(\prod_{\{
(m,j)(m^{\prime},j^{\prime})\}\in\gamma }e^{il_{m,j}\tau k_m}
e^{-il_{m,j}\tau k_{m'}}\Big)=\cr
=\Big(\prod_{\{
(m,j)(m^{\prime},j^{\prime})\}\in\gamma
}C_{l_{m,j}}a^2(\tfrac{4l_{m,j}}N)\,\delta_{l_{m,j},l_{m',j'}}\Big)\Big(\prod_{m=1}^{2p}e^{i(
l_{m,1;\gamma}+l_{m,2;\gamma }+l_{m,3;\gamma}) \tau k_m}\Big) }
$$
where $l_{m,j;\gamma}=l_{m,j}$ if
$(m,j)$ is a departing point in $\gamma$, $l_{m,j;\gamma}=-l_{m,j}$
if $(m,j)$ is an arrival point. Summing on the possible values of
$k_1,\dots,k_{2p}$, we get
$$
\displaylines{ \sum_{k_1,\dots ,k_{2p}}\prod_{\{
(m,j)(m^{\prime},j^{\prime})\}  \in\gamma
}C_{l_{m,j}}a^2(\tfrac{4l_{m,j}}N)e^{il_{m,j}\tau k_m}
e^{-il_{m,j}\tau k_{m'}}\,\delta_{l_{m,j},l_{m',j'}}=\cr
=\Big(\prod_{\{  (m,j)(m^{\prime},j^{\prime})\}\in\gamma
}C_{l_{m,j}}a^2(\tfrac{4l_{m,j}}N)\,\delta_{l_{m,j},l_{m',j'}}\Big)
\Big(\prod_{m=1}^{2p}D_N([ l_{m,1;\gamma}+l_{m,2;\gamma}+
l_{m,3;\gamma}]\tau)\Big)=\cr =\Big(\smash{\underbrace{\prod_{\{
(m,j)(m^{\prime},j^{\prime})\}\in\gamma}
\delta_{l_{m,j},l_{m',j'}}}_{=\delta(\gamma;l_{1,1},\dots,l_{2p,3})}}
\Big)\Big(\prod_{\{ (m,j)(m^{\prime},j^{\prime})\}\in\gamma
}C_{l_{m,j}}a^2(\tfrac{4l_{m,j}}N)\Big) \times\hfill\cr
\hfill\times\Big(\prod_{m=1}^{2p} D_N([
l_{m,1;\gamma}+l_{m,2;\gamma}+ l_{m,3;\gamma}]\tau)\Big)=\cr
=\delta(\gamma;l_{1,1},\dots,l_{2p,3})\prod_{m=1}^{2p}\Big(C_{l_{m,1}}
a^2(\tfrac{4l_{m,1}}N)C_{l_{m,2}}a^2(\tfrac{4l_{m,2}}N)
C_{l_{m,3}}a^2(\tfrac{4l_{m,3}}N) \Big)^{1/2}\times\hfill\cr
\hfill\times D_N([ l_{m,1;\gamma}+l_{m,2;\gamma}+
l_{m,3;\gamma}]\tau). \cr }
$$
Let us define
$$
\displaylines{
X_{l_{m,1},l_{m,2},l_{m,3};\gamma}=\cr
=\frac{1}{N^{2}\sigma_N^3}\Big(C_{l_{m,1}}a^2(\tfrac{4l_{m,1}}N)
C_{l_{m,2}}a^2(\tfrac{4l_{m,2}}N)
C_{l_{m,3}}a^2(\tfrac{4l_{m,3}}N) \Big)^{1/2}D_N([
l_{m,1;\gamma}+l_{m,2;\gamma}+ l_{m,3;\gamma}]\tau)\ .\cr }
$$
In conclusion
$$
ES_{N}^{2p}=\sum_{\gamma\in\Gamma(2p,3)}
\sum_{l_{1,1},\dots,l_{2p,3}}\delta(\gamma;l_{1,1},\dots,l_{2p,3})
\prod_{m=1}^{2p}X_{l_{m,1},l_{m,2},l_{m,3};\gamma}.
$$
%
Recall that $\delta(\gamma;l_{1,1},\dots,l_{2p,3})=0$ unless
$l_{m,j}=l_{m',j'}$ for every $(m,j)(m^{\prime},j^{\prime})\}\in\gamma$.
The proof is done by proving that the leading contribution to $ES_{N}^{2p}$
is given by paired
diagrams whereas the non paired ones are negligible in the asymptotics.
This is made explicit in the two following lemmas.
%
\begin{lemma}
For the terms corresponding to the paired diagrams $\gamma\in\Gamma_{P}(2p,3)$ we
have%
$$
\sum_{\gamma\in\Gamma_P(2p,3)}
\sum_{l_{11},\dots,l_{2p,3}}\delta(\gamma;l_{1,1},\dots,l_{2p,3})
\prod_{m=1}^{2p}X_{l_{m,1}l_{m,2}l_{m,3};\gamma}
=(2p-1)!!\sigma_{S_N}^{2p}\ .
$$
\end{lemma}
\proof Remark first that the number of possible ways of
partitioning the $2p$ rows of the diagram into subsets of cardinality $2$ is
exactly $(2p-1)!!$. Also, in every diagram $\gamma\in\Gamma_P(2p,3)$ the
contribution of two paired rows is exactly $\sigma_{S_N}^{2}$.

\hfill$\square$

\noindent To conclude the proof, we need only to show that the terms corresponding to
all remaining diagrams $\gamma\notin\Gamma_{P}(2p,3)$ are asymptotically of
smaller order.

\begin{lemma} \label{lemboun} All terms corresponding to diagrams with connected
components of order larger than 2 ($\gamma\notin\Gamma_{P}(2p,3))$
are of order $O(\frac{\log N}{\sqrt{N}})$.
\end{lemma}
\proof We focus on any two nodes that are connected but not
paired; such two nodes certainly exist, otherwise $\gamma\in\Gamma_{P}.$ We
consider the case where there is a single edge linking these two nodes; the
proof in the remaining case is entirely analogous. Without loss of
generality, we label edges and vertices in such a way that
the edge connecting these two nodes is labelled $\{(1,1),(2,1)\}$.

As in Lemma 3.1 of \cite{mari2006b}, we apply iteratively the
Cauchy-Schwartz inequality to show that%
$$
\displaylines{
\Big| \sum_{l_{11},\dots,l_{2p,3}}\delta(\gamma;l_{1,1},\dots,l_{2p,3})
\prod_{m=1}^{2p}X_{l_{m1},l_{m2},l_{m3};\gamma}\Big|  \leq\cr
\le\prod_{m=1}^{2p}\bigg(
\sum_{l_{m,1}l_{m,2}l_{m,3}}X_{l_{m,1}l_{m,2}
l_{m,3};\gamma}^2\bigg)  ^{1/2}\times\bigg(  \sum_{l_{1,2
}l_{1,3}l_{2,2}l_{2,3}}Y_{l_{1,2
}l_{1,3}l_{2,2}l_{2,3};\gamma}^2\bigg)  ^{1/2}\ ,\cr
}
$$
where
$$
\displaylines{
Y_{l_{1,2}l_{1,3}l_{2,2}l_{2,3};\gamma}=
\sum_{l_{1,1},l_{2,1}} \delta_{l_{1,1},-l_{2,1}}
X_{l_{1,1}l_{1,2}l_{1,3};\gamma}X_{l_{2,1}l_{2,2}l_{2,3};\gamma}=\cr
=\frac1{N^3\sigma_N^6}\sqrt{\,C_{l_{1,2}}a^2(\tfrac{4l_{1,2}}{N})C_{l_{1,3}}%
a^2(\tfrac{4l_{1,3}}{N})C_{l_{2,2}}a^2(\tfrac{4l_{2,2}}N)C_{l_{2,3}}a%
^{2}(\tfrac{4l_{2,3}}{N})}\,\times\hfill\cr
\hfill\times\frac1N\sum_{l_{1,1}}C_{l_{1,1}}a^2(\tfrac{4l_{1,1}}N)
D_{N}((l_{1,1;\gamma}+l_{1,2;\gamma}+l_{1,3;\gamma})\tau)D_N((-l_{1,1;\gamma}%
+l_{2,2;\gamma}+l_{2,3;\gamma})\tau)\ .\cr
}
$$
Now%
$$
\displaylines{
\sum_{l_{m,1}l_{m,2}l_{m,3}}X_{l_{m,1}l_{m,2
}l_{m,3};\gamma}^2=\cr
=\frac1{N^3\sigma_N^6}\sum_{l_{m,1}l_{m,2}l_{m,3}}
C_{l_{m,1}}a^2(\tfrac{4l_{m,1}}{N})C_{l_{m,2}}a
^2(\tfrac{4l_{m,2}}N)C_{l_{m,3}}a^2(\tfrac{4l_{m,3}%
}{N})\times\hfill\cr
\hfill\times\frac1N\, D_N([  l_{m,1;\gamma}+
l_{m,2;\gamma}+l_{m,3;\gamma}]  \tau)^2=\cr
=O\Bigg(  \frac1{N^3}\sum_{l_{m,1;\gamma}l_{m,2;\gamma}l_{m,3;\gamma}}%
K_{N}([  l_{m,1;\gamma}+l_{m,2;\gamma}+l_{m,3;\gamma}]\tau)\Bigg)=O(1)\ .\cr
}
$$
On the other hand%
$$
\displaylines{
\sum_{l_{1,2}l_{1,3}l_{2,2}l_{2,3}}Y_{l_{1,2}l_{1,3}l_{2,2}l_{2,3};\gamma}^2=\cr
=\frac{1}{N^{4}}\sum_{l_{1,2}l_{1,3}l_{2,2}l_{2,3}
}\!\!\!\!\!\!\!\!\frac{C_{l_{1,2}}a^2(\frac{4l_{1,2}}{N})C_{l_{1,3}}%
a^2(\frac{4l_{1,3}}{N})C_{l_{2,2}}a^2(\frac{4l_{2,2}}{N})C_{l_{2,3}}
a^2(\frac{4l_{2,3}}{N})}
{\sigma_{N}^{8}}\times\hfill\cr
\hfill\times\Bigg[  \frac{1}{N}\sum_{l_{1,1}}\frac{C_{l_{1,1}}{a}%
^{2}(\frac{4l_{1,1}}{N})}{\sigma_{N}^{2}}\frac{D_N((
l_{1,1;\gamma}+l_{1,2;\gamma}+l_{1,3;\gamma})\tau)D_N((-l_{1,1;\gamma}+
l_{2,2;\gamma}+l_{2,3;\gamma})\tau)}N\Bigg]  ^{2}.\cr
}
$$
Now we observe that
$$
\displaylines{
D_{N}\big((l_{1,1;\gamma}+l_{1,1;\gamma}+
l_{1,3;\gamma})\tau\big)D_N%
\big((-l_{1,1;\gamma}+l_{2,2;\gamma}+l_{2,3;\gamma})\tau\big)=\cr
=\sum_{k_{1}=1}^{N}
e^{(l_{1,1;\gamma}+l_{1,2;\gamma}+l_{1,3;\gamma})\tau k_1} \sum_{k_2=1}^{N}
e^{(-l_{1,1;\gamma}+l_{2,2;\gamma}+l_{2,3;\gamma})\tau k_2} \cr
=\sum_{k_1=1}^N
e^{(l_{1,1;\gamma}+l_{1,2;\gamma}+l_{1,3;\gamma})\tau k_1}
\sum_{u=1}^{N}e^{(l_{1,1;\gamma}+l_{2,2;\gamma}+l_{2,3;\gamma})\tau (k_1-u)} \cr
=D_N\big((l_{1,2;\gamma}+l_{1,3;\gamma}+
l_{2,2;\gamma}+l_{2,3;\gamma})\tau\big)D_{N}\big((-l_{1,1;\gamma}+
l_{2,2;\gamma}+l_{2,3;\gamma})\tau\big),\cr
}
$$
whence%
$$
\displaylines{
\Big|\frac1{N^2\sigma_N^2} \sum_{l_{1,1}}C_{l_{1,1}}a^2
(\tfrac{4l_{1,1}}N)D_N((l_{1,1;\gamma}+l_{1,2;\gamma}+l_{1,3;\gamma})\tau)
D_N((-l_{1,1;\gamma}+l_{2,2;\gamma}+l_{2,3;\gamma})\tau)\Big| \le\cr
\leq\Big|\frac 1N
D_N((l_{1,2;\gamma}+l_{1,3;\gamma}+l_{2,2;\gamma}+l_{2,3;\gamma})\tau)\Big|\times\hfill\cr
\hfill\times  \Big|  \frac1N\sum_{l_{1,1}}%
\frac1{\sigma_N^2}C_{l_{11;\gamma}}a^2(\tfrac{4l_{1,1}}N)
D_N((-l_{1,1;\gamma}+l_{2,2;\gamma}+l_{2,3;\gamma})\tau)\Big|\le
\cr
\leq\frac{C}{\left|  l_{1,2;\gamma}+l_{1,3;\gamma}+l_{2,2;\gamma}%
+l_{2,3;\gamma}\right|  +1}
  \frac1N\sum_{l_{1,1;\gamma}}\frac
{N}{|-l_{1,1;\gamma}+l_{2,2;\gamma}+l_{2,3;\gamma}|+1}\le \cr
\leq C\,\frac{\log N}{\left|  l_{1,2;\gamma}+l_{1,3;\gamma}+
l_{2,2;\gamma}+l_{2,3;\gamma}\right|+1}\cdotp\cr
}
$$
Thus we can conclude that%
$$
\displaylines{
\bigg(  \sum_{l_{1,2}l_{1,3}l_{2,2}l_{2,3}}
Y_{l_{1,2}l_{1,3}l_{2,2}l_{2,3};\gamma}^2
\bigg)^{1/2}\le\cr
\leq c\Bigg(\frac{1}{N^{4}}\sum_{l_{1,2}l_{1,3}l_{2,2}l_{2,3}}
\frac{\log^{2}N}{\big(  |  l_{1,2;\gamma}+l_{1,3;\gamma}
+l_{2,2;\gamma}+l_{2,3;\gamma}|  +1\big)^2}\Bigg)^{1/2}
\leq c\,\frac{\log N}{\sqrt{N}}\cdotp\cr
}
$$
\qquad\hfill$\square$
%
%
\section{Studentized statistics}\label{stud}

\subsection{Estimation of $\sigma_{N}^{2}$.}

The statistics described in the previous sections can be impossible to
compute in
practice, as they depend on the correlation
structure of the field, which is in general unknown. In this Section, we show
how asymptotic variances can be consistently estimated from the data, in the
presence of observations at higher and higher resolution. We start from the
variance of the wavelets coefficients, which we recall is given by
\[
\sigma_{N}^{2}:=\frac2N\sum_{N/8\leq l\leq N/2}C_{l}a^2(\tfrac{4l}%
{N})\text{ .}%
\]
An obvious estimator is provided by%
\[
\widehat{\sigma}_{N}^{2}:=\frac2N\sum_{N/8\leq l\leq N/2}
|w_l|^2%
a^2(\tfrac{4l}{N})\ .%
\]
Of course, $\widehat{\sigma}_{N}^{2}$ is unbiased and mean square consistent
for $\sigma_{N}^{2},$ in the trivial sense that both converge to zero as $N$
diverges. The following result is stronger.

\begin{lemma}\label{conv-stud-var1}
Under assumptions A1 and B, as $N\rightarrow\infty$, we have%
\[
\lim_{N\to\infty} \frac{\widehat{\sigma}_N^2}{\sigma_N^2}=1
\qquad \mbox{in }L^2\ .
\]
\end{lemma}
\proof It is immediate to see that%
\[
E\Big(\frac{\widehat{\sigma}_N^2}{\sigma_N^2}\Big)=E\bigg(  \frac{\sum_{N/8\leq l\leq N/2}|w_{l}|^{2}a^2(\frac{4l}N)}
{\sum_{N/8\leq l\leq N/2}C_{l}a^2(\frac{4l}N)}\bigg)
\equiv1
\]
and%
\begin{align*}
\Var(\widehat{\sigma}_N^2)&=\sum_{N/8\leq l\leq N/2}\Var(|w_l|^2)a^2(\tfrac{4l}%
{N})=2\sum_{N/8\leq l\leq N/2}C_l^2a^4(\tfrac{4l_1}N)\ .
\end{align*}
It is sufficient now to note that, under Assumption A1,
\[
0<c_{1}\leq\frac{C_{l}}{C_{N/4}}\leq c_{2}<\infty,\qquad\text{for }%
N=2^{j+2}\text{ and for all}\ l\in\lbrack N/8,N/2])\text{ ,}%
\]
whence%
\[
\Var\Big(\frac{\widehat{\sigma}_N^2}{\sigma_N^2}\Big)=\frac{\sum_{l}C_{l}^{2}{a}^{4}(\frac{4l_1}{N})}{\left\{  \sum_{l}C_{l}%
a^2(\frac{4l}{N})\right\}^2}=O\left(
\frac{\sum_la^4(\frac {4l}N)}{\left\{
\sum_{l}a^2(\frac{4l}{N})\right\}^2}\right)\ .
\]
Remark now that $\sum_la^4(\frac {4l}N)\sim N\int_{1/2}^1a^4(t)\, dt$,
$\sum_la^2(\frac {4l}N)\sim N\int_{1/2}^1a^2(t)\, dt$, so that the left-hand
term tends to $0$ as $\frac 1N$.

\qquad\hfill$\square$

\subsection{Estimation of the variance for Skewness and kurtosis}

We now go on with the estimation for the sample variance for the
statistics $S_{N}$ and $U_{N}.$ We note first that under Gaussianity
$( \left| w_{l}\right|  ^{2}/C_{l})_l  $ is a sequence of
independent and identically distributed exponential random variables
with mean unity; we
define for all $p\in\mathbb{N}$%
\[
\delta_{l_1l_2\dots l_p}:=E\bigg( \prod_{l=l_1}^{l_{p}}\frac{\left|
w_{l}\right|  ^{2}}{C_{l}}\bigg)  \text{ ;}%
\]
to be quite explicit we have for instance%
\begin{align*}
\delta_{l_1l_2}  &  =\begin{cases}
1&\text{ if }|l_1|\neq |l_2|\\
2&\text{ if }|l_1|=|l_2|%
\end{cases}
\\
\delta_{l_1l_2l_3}  &  =\begin{cases}
1&\text{ if }|l_1|, |l_2|, |l_3| \text{ are distinct }\\
2&\text{ if among }|l_1|,|l_2|,|l_3|\text{ two are equal
and the third is different }\\
6&\text{ if }|l_1|=|l_2|=|l_3|%
\end{cases}
\end{align*}
In view of (\ref{varsk})--(\ref{varkur2}), a natural proposal is to
consider
for the Skewness%
\begin{equation}
\widehat{\sigma}_{S_{N}}^{2}:=
\frac{12\pi}{N^{3}\widehat{\sigma}_{N}^{6}}%
\sum_{l_1l_2l_3}\frac{1}{\delta_{l_1l_2l_3}}
|w_{l_1}|^{2}a^2(\tfrac{4l_1}{N}) |w_{l_2}|^{2}a^2(\tfrac{4l_2}{N})
|w_{l_3}|^{2}a^2(\tfrac{4l_3}{N})
K_{N}(\tfrac{2\pi}N(l_1+l_2+l_3))\label{estskew}%
\end{equation}
and for the Kurtosis%
\begin{align}
\displaystyle\widehat{\sigma}_{U_{N}}^{2}&:=\widehat{\sigma}_{1U_{N}}^{2}+\widehat{\sigma
}_{2U_{N}}^{2}\text{ ,}\label{estkur0} \\
\displaystyle\widehat{\sigma}_{1U_{N}}^{2}&:=\frac{72}{\widehat\sigma_{N}^{4}}
\frac{2\pi}{N^{2}}
\sum_{l_1l_2}\frac{1}{\delta_{l_1l_2}}|w_{l_1}|^{2}
a^2(\tfrac{4l_1}{N})|w_{l_2}|^{2}a^2(\tfrac{4l_2}{N})K_{N}(\tfrac
{2\pi}N(l_1+l_2))\text{ ,}\label{estkur1}\\
\displaystyle\widehat{\sigma}_{2U_{N}}^{2}&:=\frac{24}{\widehat\sigma_{N}^{8}}\frac{2\pi}{N^{4}}%
\sum_{l_1l_2l_3l_{4}}\frac{1}{\delta_{l_1l_2l_3l_{4}}}\bigg\{
\prod_{l=l_1}^{l_{4}}|w_{l}|^{2}a^2(\tfrac{4l}{N})\bigg\}
K_{N}(\tfrac{2\pi}N(l_1+l_2+l_3+l_{4})).\label{estkur2}
\end{align}

\begin{rem}Using the properties of Fej\'{e}r's kernel
recalled in \S \ref{sk-ku-ssec}, in the summations above most
terms vanish. From a computational point of view more tractable
expressions can be derived in the spirit of \S \ref{sk-ku-ssec}.
In particular it holds
\begin{equation}\label{estskew-due}
\widehat\sigma_{1U_{N}}^{2}:=\frac{72}{\widehat\sigma_{N}^{4}}\frac{1}{N}
\sum_{l}|w_{l}|^4a^4(\tfrac{4l}{N})\ .%
\end{equation}
\end{rem}
\begin{lemma}\label{conv-stud-stat-var}
Under Assumptions A1 and B, as $N\rightarrow\infty,$ we have%
\[
\frac{\widehat{\sigma}_{S_N^2}}{\sigma_{S_N^2}}\enspace
\mathop{\to}_{N\to\infty}^P\enspace 1,\qquad
\frac{\widehat{\sigma}_{1U_N^2}}{\sigma_{1U_N}^2}\enspace
\mathop{\to}_{N\to\infty}^P\enspace 1, \qquad
\frac{\widehat{\sigma}_{2U_N^2}}{\sigma_{2U_N}^2}\enspace
\mathop{\to}_{N\to\infty}^P\enspace 1\ .%
\]
\end{lemma}
\proof We give the proof for $\widehat{\sigma}_{1U_{N}}%
^{2},\widehat{\sigma}_{2U_{N}}^{2}$ only, as the remaining case is
entirely analogous (indeed slightly simpler). Let us denote
$\widetilde{\sigma}_{1U_N}^2 =\widehat{\sigma}_{1U_N}^2\cdot
\widehat{\sigma}_N^4/\sigma_N^4$, that is the same as in
(\ref{estskew}) with $\widehat{\sigma}_N$ replaced by ${\sigma}_N$.
As $E\Big(  \frac{1}{\delta_{l_1l_2}}|w_{l_1}|^{2}|w_{l_2}|^{2}
\Big) =C_{l_1}C_{l_2}$ for every $l_1,l_2$, it is clear that
\begin{equation}\label{kurvar12}
E\bigg(\frac{\widetilde{\sigma}_{1U_N}^2} {\sigma_{1U_N}^2}\bigg)
=1 .%
\end{equation}
Moreover, by the alternate expression (\ref{estskew-due}) and in
view of Remark \ref{conv-var-coeff},
$$
\Var(\widetilde{\sigma}_{1U_N}^2)=\frac{72^2}{\sigma_N^8}\frac1{N^2}
\sum_{l}\Var(|w_{l}|^4)a^8(\tfrac{4l}N)= \frac
{c_0}N\,\frac{C_{N/4}^4} {\sigma_N^8}\, \frac
1N\sum_{l}a^8(\tfrac{4l}N)\sim \frac{c_1}N\cdotp
$$
%
As we know that under Assumption A1 $\sigma_{1U_N}^2$ is bounded
away from zero (see Remark \ref{conv-var-skku}), this implies that
$\Var(\widetilde{\sigma}_{1U_N}^2/\sigma_{1U_N}^2)\to 0$ as
$N\to\infty$.

The argument for $\widehat{\sigma}_{2U_{N}}^{2}$ is similar; indeed
if we define $\widetilde{\sigma}_{2U_N}^2$ in analogy with
$\widetilde{\sigma}_{1U_N}^2$, it is immediate that
$E[\widehat{\sigma}_{2U_{N}}^{2}/\sigma_{2U_N}^{2}]=1$. On the other
hand, note that the summands in $\widetilde{\sigma}_{2U_N}^2$ have a
martingale-difference structure on the lattice $\mathbb{Z}^{3}$ (see
Poghosyan and Roelly \cite{MR1629903} e.g.), whence, in view of
Assumption A
$$
\displaylines{
\Var\bigg(  \frac{\widehat{\sigma}_{2U_{N}}^{2}}{\sigma_{2U_{N}}^{2}%
}\bigg)=\cr =O\bigg( \Var\bigg\{ \frac{24}{\sigma_{N}^{8}}\frac{2\pi
}{N^4}\sum_{l_1l_2l_3l_{4}}\frac{1}{\delta_{l_1l_2l_3l_{4}}%
}\bigg[  \prod_{l=l_1}^{l_{4}}|w_{l}|^{2}a^2(\tfrac{4l}{N}) \bigg]
K_{N}(\tfrac{2\pi}N(l_1+l_2+l_3+l_{4}))\bigg\} \bigg) \cr
  =O\bigg(  \frac{1}{N^{8}}\sum_{l_1l_2l_3l_{4}}K_{N}^{2}
  (\tfrac{2\pi}N(l_1+l_2+l_3+l_{4}))\Var\bigg\{
\prod_{l=l_1}^{l_{4}}%
\frac{|w_{l}|^{2}}{\sigma_{N}^{2}}\bigg\}  \bigg) \cr
  =O\bigg(  \frac{1}{N^{6}}\sum_{l_1l_2l_3l_{4}}\Var\bigg\{
\prod_{l=l_1}^{l_{4}}\frac{|w_{l}|^{2}}{\sigma_{N}^{2}}\bigg\}
\bigg)
=O(\tfrac{1}{N^{2}})=o(1)\text{ .}\cr%
}
$$
We have thus proved that
$\widetilde{\sigma}_{1U_N}^2/\sigma_{1U_N}^2\to 1$ and
$\widetilde{\sigma}_{2U_N}^2/\sigma_{2U_N}^2\to 1$ as $N\to\infty$
in $L^2$. Therefore, in view of Lemma \ref{conv-stud-var1},
$\widetilde{\sigma}_{1U_N}^2/\sigma_{1U_N}^2\to 1$ and
$\widetilde{\sigma}_{2U_N}^2/\sigma_{2U_N}^2\to 1$ as $N\to\infty$
in probability. The rest of the proof is quite similar.

\hfill$\square$
\newline\noindent%
%
As an immediate consequence of Theorem \ref{TLC} and Lemma
\ref{conv-stud-stat-var} we have the following.
\begin{theorem}
\label{TLC2}Under the assumptions A1 and B, as $N\rightarrow\infty$%
\[
\begin{pmatrix}
\frac1{\widehat{\sigma}_{S_N}}S_N\\
\frac1{\widehat{\sigma}_{U_N}}U_N%
\end{pmatrix}
\enspace\mathop{\rightarrow}^{\mathcal D}_{N\to\infty}
\enspace N(0,I_{2})\text{ .}%
\]
\end{theorem}

\section{Aliasing}\label{alias}

The tests provided in the sections above are based on the wavelet
coefficients
\[
\beta_{Nk}:=\frac1{2\pi}\int_{-\pi}^{\pi}X(\vartheta)\psi_{Nk}
(\vartheta)d\vartheta,\hbox to 0pt{\qquad \hfill$k=0,1,\dots,N-1$}\ .
\]
In practice, $\beta_{Nk}$ will be approximated by the sums
\[
\widetilde{\beta}_{Nk}=\frac1M\sum_{m=0}^{M-1}X(\tfrac{2\pi
}M\,m)\psi _{Nk}(\tfrac{2\pi }M\,m)\ .
\]
The purpose of this section is to prove that, if our data have
enough high frequencies, this approximation does not affect the
test.

We need to strengthen our previous assumptions as follows.
\smallskip

\noindent\textbf{Assumption C} As $l\rightarrow\infty$ we have%
\[
C_{l}=L(l)l^{-\alpha},\text{ }\alpha>1
\]
where $L(l)$ denotes a slowly varying function \cite{MR1015093},
which we assume to be bounded and bounded away from zero.
\smallskip

\noindent\textbf{Assumption D }$N$ is such that
\[
\frac{1}{N}+\frac{N}{M^{\alpha/(\alpha+1)}}\rightarrow0\text{ as }%
M\rightarrow\infty\text{ .}%
\]
\smallskip

\noindent Assumption C is mild, entailing simply a regular behaviour
of the angular power spectrum at infinity. Assumption D is a sort of
bandwidth condition, suggesting that the frequencies that we can use
fruitfully for statistical inference must grow more slowly than the
sampling rate. The condition become less and less tight the faster
the angular correlation function decays to zero: for instance if
$\alpha=4$ we must impose $N=o(M^{4/5})$. In practice $\alpha$ can
be estimated from the data; the most cautious choice can be
$M=o(\sqrt{N}),$ as $\alpha>1$ is implied by the finite variance of
the field.

We have the following result.
\begin{proposition}
Under assumptions A1,B,C and D the result of Theorem \ref{TLC} remains true
when replacing the $\beta_{Nk}$'s by the $\widetilde{\beta}_{Nk}$'s.
\end{proposition}
\proof It holds,
\[
\beta_{Nk}-\widetilde{\beta}_{Nk}=\sum_{|\ell|>M-\frac{N}{2}}
w_{\ell}\Big(\frac
{1}{M}\sum_{m=0}^{M-1}\psi_{Nk}(\tfrac{2\pi m}{M})e^{i2\pi m\frac{\ell}{M}%
}\Big)\text{ .}%
\]
Therefore,
\begin{align*}
E\big|  \beta_{Nk}-\widetilde{\beta}_{Nk}\big|  ^{2}  &  =E\bigg(
\sum_{|\ell|>M-\frac{N}{2}}|w_{\ell}|^{2}\frac{1}{M^{2}}\bigg|  \sum
_{m=0}^{M-1}\psi_{Nk}(\frac{2\pi m}{M})e^{i2\pi m\frac{\ell}{M}}\bigg|
^{2}\bigg)
=\sum_{|\ell|>M-\frac{N}{2}}C_{l}\left|  A_{M(l)}\right|  ^{2},
\end{align*}
where
\[
A_{M(l)}^{2}=\frac{1}{M^{2}}\Big|\sum_{m=0}^{M-1}\psi_{Nk}(\tfrac{2\pi m}%
{M})e^{i2\pi m\frac{l}{M}}\Big|^{2}\leq\frac{c_{k}}{|1+[l]_{M}|^{k}}\text{ for all
}k>0\text{ , some }c_{k}>0\text{ ,}%
\]
in view of Theorem (\ref{concen}), which implies
\[
\Big|  \sum_{m=0}^{M-1}\psi_{Nk}(\tfrac{2\pi m}{M})e^{i2\pi m\frac{l}{M}%
}\Big|  \leq\frac{c_{k}M}{|1+M[l/M]_{2\pi}|^{k}}\text{ .}%
\]
Under Assumptions C and D we have easily
$$
\displaylines{
\sum_{|\ell|>M-\frac{N}{2}}C_{l}\left|  A_{M(l)}\right|  ^{2}=O\bigg(
C_{M}\sum_{|\ell|>M-\frac{N}{2}}\frac{C_{l}}{C_{M}}\left|  A_{M(l)}\right|
^{2}\bigg) \cr
=O\bigg(  C_{M}\sum_{u=1}^{\infty}\sum_{v=-M}^{M}\frac{C_{uM+v}}{C_{M}%
}\left|  A_{M(uM+v)}\right|  ^{2}\bigg)=\cr
=O\bigg(  M^{-\alpha}\sum_{u=1}^{\infty}u^{-\alpha}\sum_{v=-M}^{M}\frac
{1}{(1+|v|)^{k}}\bigg)  =O\left(  M^{-\alpha}\right)\cr
}
$$
and
\begin{equation}
\frac{E(|\beta_{Nk}-\widetilde{\beta}_{Nk}|^{2})}{E|\beta_{Nk}|^{2}}\leq
c\frac{M}{\sigma_{N}^{2}}^{-\alpha}\leq c\Big(  \frac{M}{N}\Big)
^{-\alpha}. \label{bound1}%
\end{equation}
The result then follows from the following lemma.
\begin{lemma}
\label{alias-lem}For $X_{in}$ and $Y_{in}$ mutually Gaussian centered random
variables, set
\[
c_{n}:=\frac{E[(X_{in}-Y_{in})^{2}]}{E[X_{in}^{2}]}.
\]
Let us assume that the sequence
\begin{equation}
\frac{\sum_{i=1}^{n}X_{in}^{3}}{\sqrt{\Var\Big(  \sum_{i=1}^{n}X_{in}%
^{3}\Big)  }} \label{conv}%
\end{equation}
converges in distribution to a variable $X$, where for $\gamma_{1},\gamma
_{2}>0$
\begin{equation}
\gamma_{1}\leq\sqrt{\frac{1}{n}\Var\bigg\{  \sum_{\smash{i=1}}^{n}X_{in}^{3}\bigg\}
}\leq\gamma_{2}\text{ and }c_{n}={o}(\frac{1}{n})\text{ }. \label{concond}%
\end{equation}
Then
\begin{equation}
\frac{\sum_{i=1}^{n}Y_{in}^{3}}{\sqrt{\Var\left\{  \sum_{i=1}^{n}X_{in}%
^{3}\right\}  }} \label{resconv}%
\end{equation}
also converges in distribution to $X$. The same result is true if we replace
in (\ref{conv}),(\ref{concond}) $X_{in}^{3}$ by $X_{in}^{4}-EX_{in}^{4}$ and
replace in (\ref{resconv}), $Y_{in}^{3}$ by $Y_{in}^{4}-EY_{in}^{4}$.
\end{lemma}
\proof We shall actually prove a stronger result, namely
\[
\lim_{n\rightarrow\infty}E\left|  \frac{\sum_{i=1}^{n}X_{in}^{3}}%
{\sqrt{\Var\left\{  \sum_{i=1}^{n}X_{in}^{3}\right\}  }}-\frac{\sum_{i=1}%
^{n}Y_{in}^{3}}{\sqrt{\Var\left\{  \sum_{i=1}^{n}X_{in}^{3}\right\}  }}\right|
=0\text{ .}%
\]
Using $(x-y)^{3}=x^{3}-y^{3}-3x(x-y)^{2}+3x^{2}(x-y)$, we get :
$$
\displaylines{
E|X_{in}^{3}-Y_{in}^{3}|    \leq E|X_{in}-Y_{in}|^{3}+3E|X_{in}-Y_{in}%
|^{2}|X_{in}|+3E|X_{in}-Y_{in}||X_{in}|^{2}\le\cr
 \leq E|X_{in}-Y_{in}|^{3}+3[E|X_{in}-Y_{in}|^{4}]^{\frac{1}{2}}%
[E|X_{in}|^{2}]^{\frac{1}{2}}+3[E(X_{in}-Y_{in})^{2}]^{\frac{1}{2}}%
[E|X_{in}|^{4}]^{\frac{1}{2}}.\cr
}
$$
Now, when $Z$ is a Gaussian random variable, for $h\;\in\mathbb{N}^{\ast}$,
\[
E|Z|^{h}=\sigma_{h}[E|Z|^{2}]^{\frac{h}{2}},
\]
where $\sigma_{h}$ is the $h$-moment of the standard gaussian distribution,
centered and with variance 1. Therefore, we have $E|X_{in}-Y_{in}|^{k}%
\leq\sigma_{h}c_{n}^{\frac{h}{2}}[EX_{in}^{2}]^{\frac{h}{2}}$, and
$$
\displaylines{
\sum_{i=1}^{n}E|X_{in}^{3}-Y_{in}^{3}|    \leq\sigma_{3}c_{n}^{\frac{3}{2}%
}\sum_{i=1}^{n}[EX_{in}^{2}]^{\frac{3}{2}}+3\sigma_{4}^{\frac{1}{2}}c_{n}%
[\sum_{i=1}^{n}[EX_{in}^{2}]^{\frac{3}{2}}+3c_{n}^{\frac{1}{2}}\sigma
_{4}^{\frac{1}{2}}\sum_{i=1}^{n}[EX_{in}^{2}]^{\frac{3}{2}}\cr
\leq\sum_{i=1}^{n}[EX_{in}^{2}]^{\frac{3}{2}}\{\sigma_{3}c_{n}^{\frac{3}%
{2}}+3\sigma_{4}^{\frac{1}{2}}c_{n}^{\frac{1}{2}}+3\sigma_{4}^{\frac{1}{2}%
}c_{n}\}.\cr
}$$
Therefore
$$
\frac{\sum_{i=1}^{n}E|X_{in}^{3}-Y_{in}^{3}|}{\sqrt{\Var\left\{  \sum_{i=1}%
^{n}X_{in}^{3}\right\}  }} \leq C\frac{\sum_{i=1}^{n}[EX_{in}^{2}%
]^{\frac{3}{2}}c_{n}^{\frac{1}{2}}}{\sqrt{n}}
\leq C\sqrt{nc_{n}}=o(1),\qquad\text{as }n\rightarrow\infty.%
$$
This proves the result for the third power. As for the forth one, we proceed
in the same way, and prove using the same path,
\begin{align*}
\displaystyle &|x^{4}-y^{4}|   \leq4|x|^{3}|x-y|+6x^{2}|x-y|^{2}+4|x]|x-y|^{3}+|x-y|^{4}\\
\displaystyle &E|X_{in}^{4}-Y_{in}^{4}|    \leq C\sqrt{c_{n}}[EX_{in}^{2}]^{2}\\
\displaystyle&\frac{\sum_{i=1}^{n}E|X_{in}^{4}-Y_{in}^{4}|}{\sqrt{\Var\left\{  \sum_{i=1}%
^{n}X_{in}^{4}\right\}  }}    \leq C\sqrt{c_{n}n}=o(1)\text{ .}%
\end{align*}
\hfill$\square$

\noindent Our final result in this Section extends the analysis of the aliasing effect
to studentized statistics. In particular, in the previous section the
variances of Skewness and Kurtosis where estimated on the basis of the
spectral coefficients $\left\{  w_{l}\right\}  ,$ which are obtained as
\[
w_{l}=\frac 1{2\pi}\int_{-\pi}^{\pi}X(\vartheta)e^{-il\vartheta}\,d\vartheta\text{ },\text{
}l=1,2,\dots,N\ .
\]
As before, in practice these Fourier coefficients will be approximated by
interpolations sums, i.e. for $M\geq N$ we have to consider
\[
\widetilde{w}_{l}=\frac1{M}\sum_{m=0}^{M-1}X(\tfrac{2\pi m}{M}%
)e^{-i\tfrac{2\pi m}{M}l}\text{ }.
\]
We note that%
\begin{equation}\label{alia}
\begin{array}{c}
\displaystyle\vrule height0pt depth18pt width0pt
\widetilde{w}_{l}  =\frac1{M}\sum_{m=0}^{M-1}\sum_{k=1}^{\infty}%
w_{k}e^{\frac{2i\pi }{M}mk}e^{-\frac{2i\pi }{M}ml}
=\frac1{M}\sum_{m=0}^{M-1}\sum_{k=1}^{\infty}w_{k}
e^{\frac{2i\pi}{M}(k-l)m}\\
\displaystyle=\frac1{M}\sum_{k=1}^{\infty}w_{k}D_{M}(\tfrac{2\pi}{M}(k-l)%
)=w_{l}+\sum_{k=1}^{\infty}w_{l+kM}.
\end{array}
\end{equation}
It follows, using Assumption D, that%
\begin{equation}
\frac{E\left|  w_{l}-\widetilde{w}_{l}\right|  ^{2}}{C_{l}}=\frac{1}{C_{l}%
}\sum_{k=1}^{\infty}C_{l+kM}\leq c\frac{M^{-\alpha}}{C_{l}}\sum_{k=1}^{\infty
}k^{-\alpha}\leq c\Big(  \frac{M}{N}\Big)  ^{-\alpha}=o(\tfrac1N). \label{bound2}%
\end{equation}
\begin{proposition}
Under Assumptions A1,B, C and D the result of Theorem \ref{TLC2} remains
true when replacing the $w_{l}$'s by the $\widetilde{w}_{l}$'s.
\end{proposition}
\proof It is clearly enough to show that
\begin{equation}\label{convp}
\lim_{N\to\infty}\frac{\widetilde{\sigma}_{S_{N}}^{2}-\widehat{\sigma}_{S_{N}}^{2}}%
{\widehat{\sigma}_{S_{N}}^{2}}=0\qquad\qquad
\lim_{N\to\infty}\frac{\widetilde{\sigma}_{U_{N}}^{2}%
-\widehat{\sigma}_{U_{N}}^{2}}{\widehat{\sigma}_{U_{N}}^{2}}=0
\end{equation}
in probability, where $\widehat{\sigma}_{S_{N}}^{2},\widehat{\sigma}_{U_{N}}^{2}$ are defined
in (\ref{estskew}) -(\ref{estkur0}),
\[
\widetilde{\sigma}_{S_{N}}^{2}:=\frac{12\pi}{N^{3}\widehat{\sigma}_{N}^{6}%
}\sum_{l_1l_2l_3}\frac{1}{\delta_{l_1l_2l_3}}|\widetilde{w}%
_{l_1}|^{2}a^2(\tfrac{4l_1}{N})|\widetilde{w}_{l_2}|^{2}a^2%
(\tfrac{4l_2}{N})|\widetilde{w}_{l_3}|^{2}a^2(\tfrac{4l_3}{N}%
)K_{N}(\tfrac{2\pi}N(l_1+l_2+l_3))
\]
and $\widetilde{\sigma}_{U_{N}}^{2}=\widetilde{\sigma}_{1U_{N}}^{2}%
+\widetilde{\sigma}_{2U_{N}}^{2},$ where%
\begin{align*}
\widetilde{\sigma}_{1U_{N}}^{2}&:=\frac{72}{\sigma_{N}^{4}}\frac{2\pi}{N^{2}%
}\sum_{l_1l_2}\frac{1}{\delta_{l_1l_2}}|\widetilde{w}_{l_1}|^{2}%
a^2(\tfrac{4l_1}{N})|\widetilde{w}_{l_2}|^{2}a^2(\frac{4l_2}%
{N})K_{N}(\tfrac{2\pi}N(l_1+l_2))\text{ ,}\\
\widetilde{\sigma}_{2U_{N}}^{2}&:=\frac{24}{\sigma_{N}^{8}}\frac{2\pi}{N^{4}%
}\sum_{l_1l_2l_3l_{4}}\frac{1}{\delta_{l_1l_2l_3l_{4}}}\bigg\{
\prod_{l=l_1}^{l_{4}}|\widetilde{w}_{l}|^{2}a^2(\tfrac{4l}{N})
\bigg\}K_{N}(\tfrac{2\pi}N(l_1+l_2+l_3+l_{4}))\text{ .}
\end{align*}
We focus on $\widetilde{\sigma}_{2U_{N}}^{2}$ as the other cases are strictly
analogous, indeed slightly simpler. We have
\[
\frac{1}{\widehat{\sigma}_{2U_{N}}^{2}}=O_{p}(1)
\]
and
$$
\displaylines{
E[ \widetilde{\sigma}_{2U_{N}}^{2}-\widehat{\sigma}_{2U_{N}}%
^{2}]\le\cr
\leq\frac{24}{\sigma_{N}^{8}}\frac{2\pi}{N^{4}}\sum_{l_1l_2l_3l_{4}%
}\frac{1}{\delta_{l_1l_2l_3l_{4}}}K_N(\tfrac{2\pi}{N}(l_1+l_2%
+l_3+l_{4}))
E\Big|  \prod_{l=l_1}^{l_{4}}|\widetilde{w}%
_l|^2a^2(\tfrac{4l}{N}) -\prod_{l=l_1}^{l_{4}%
}|\widetilde{w}_l|^{2}a^2(\tfrac{4l}{N})  \Big|\le \cr
\leq C\Big[  \max_{N/8\leq l_1,l_2,l_3,l_{4}\leq N/2}
\frac{1}{\sigma_{N}^{8}}E\Big|  \Big\{  \prod_{l=l_1}^{l_{4}}%
|\widetilde{w}_{l}|^{2}a^2(\tfrac{4l}{N})\Big\}  -\Big\{  \prod
_{l=l_1}^{l_4}|\widetilde{w}_l|^2a^2(\tfrac{4l}{N})\Big\}
\Big|\Big]  \text{ .}\cr%
}
$$
Now notice that%
$$
x_1x_2x_3x_4-y_1y_2y_3y_4=x_1x_2x_3(x_4-y_4)+x_1x_2(x_3-y_3)y_4
+x_1(x_2-y_2)y_3y_4+(x_1-y_1)y_2y_3y_4
$$
whence%
$$
\displaylines{ \max_{N/8\leq l_1,l_2,l_3,l_{4}\leq N/2}\bigg\{
\frac{1}{\sigma
_{N}^{8}}E\bigg|   \prod_{l=l_1}^{l_{4}}|\widetilde{w}_{l}|^{2}%
a^2(\tfrac{4l}{N})  -  \prod_{l=l_1}^{l_{4}}|w_{l}|^{2}%
a^2(\tfrac{4l}{N})  \bigg|  \bigg\}\le \cr \leq\max_{N/8\leq
l_1,l_2,l_3,l_{4}\leq N/2}\bigg\{  \bigg(
\prod_{l=l_1}^{l_{4}}a^2(\tfrac{4l}{N})\bigg) E\bigg[
\frac{ |\widetilde{w}_{l_{4}}|^{2}-|w_{l_{4}}|^{2}  }%
{\sigma_{N}^{2}}  \prod_{l=l_1}^{l_3}
\frac{|\widetilde{w}_{l}|^{2}} {\sigma_N^2}\bigg] \bigg\}  +\cr
+\max_{N/8\leq l_1,l_2,l_3,l_{4}\leq N/2}\bigg\{  \bigg(
\prod_{l=l_1}^{l_{4}}a^2(\tfrac{4l}{N})\bigg)  E\bigg|
\frac{ |\widetilde{w}_{l_3}|^{2}-|w_{l_3}|^{2}}%
{\sigma_{N}^{2}}\frac{|w_{l_{4}}|^{2}}{\sigma_{N}^{2}}
\prod_{l=l_1}^{l_2}\frac{|\widetilde{w}_{l}|^{2}}{\sigma_N^2}
\bigg|  \bigg\} +\cr
+\max_{N/8\leq l_1,l_2,l_3,l_{4}\leq N/2}\bigg\{  \bigg(
\prod_{l=l_1}^{l_{4}}a^2(\tfrac{4l}{N})\bigg)  E\bigg|  \frac
{|\widetilde{w}_{l_1}|^{2}}{\sigma_{N}^{2}}\frac{\left\{
|\widetilde
{w}_{l_2}|^{2}-|w_{l_2}|^{2}\right\}  }{\sigma_{N}^{2}}\frac{|w_{l_3%
}|^{2}|w_{l_{4}}|^{2}}{\sigma_{N}^{4}}\bigg|  \bigg\} +\cr
+\max_{N/8\leq l_1,l_2,l_3,l_{4}\leq N/2}\bigg\{  \bigg(
\prod_{l=l_1}^{l_{4}}a^2(\tfrac{4l}{N})\bigg)  E\bigg|
\frac{\left\{
|\widetilde{w}_{l_1}|^{2}-|w_{l_1}|^{2}\right\}  }{\sigma_{N}^{2}}%
\frac{|w_{l_2}|^{2}|w_{l_3}|^{2}|w_{l_{4}}|^{2}}{\sigma_{N}^{6}}
\bigg| \bigg\}\le \cr \leq \frac C{\sigma_{N}^{6}}\max_{N/8\leq
l_1,l_2,l_3,l_{4}\leq N/2}\left\{
E|\widetilde{w}_{l}|^{2}+E|w_{l}|^{2}\right\}  ^{3}%
\times\frac1{\sigma_N^2} \max_{N/8\leq l_1,l_2,l_3,l_{4}\leq N/2}
E\left| |\widetilde{w}_{l}|^{2}-|w_{l}|^{2}\right| \le  \cr \leq
C\Big(  \frac{M}{N}\Big)  ^{-\alpha}=o(1)\ .\cr }
$$
\hfill$\square$
\begin{rem}It is evident from the proof that, in order to establish
(\ref{convp}), it is sufficient to impose the minimal bandwidth condition
\[
\frac{1}{N}+\frac{N}{M}\rightarrow0\text{ as }M\rightarrow\infty\text{ .}%
\]
which is weaker than Assumption D. This is intuitively due to peculiar form that the aliasing effect assumes
for Fourier transforms in the discrete case, see (\ref{alia}).
\end{rem}
\section{Monte Carlo evidence}\label{montec}
\begin{figure}[h]
\centering
\includegraphics[height=5cm,width=9cm]{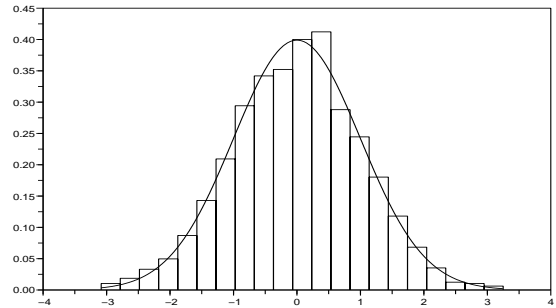}
\caption{Histogram of the skewness statistics over $1600$ simulated
fields. Here $N=2^{12}$.\label{fig-sk}}
\end{figure}

\begin{figure}[h]
\centering
\includegraphics[height=5cm,width=9cm]{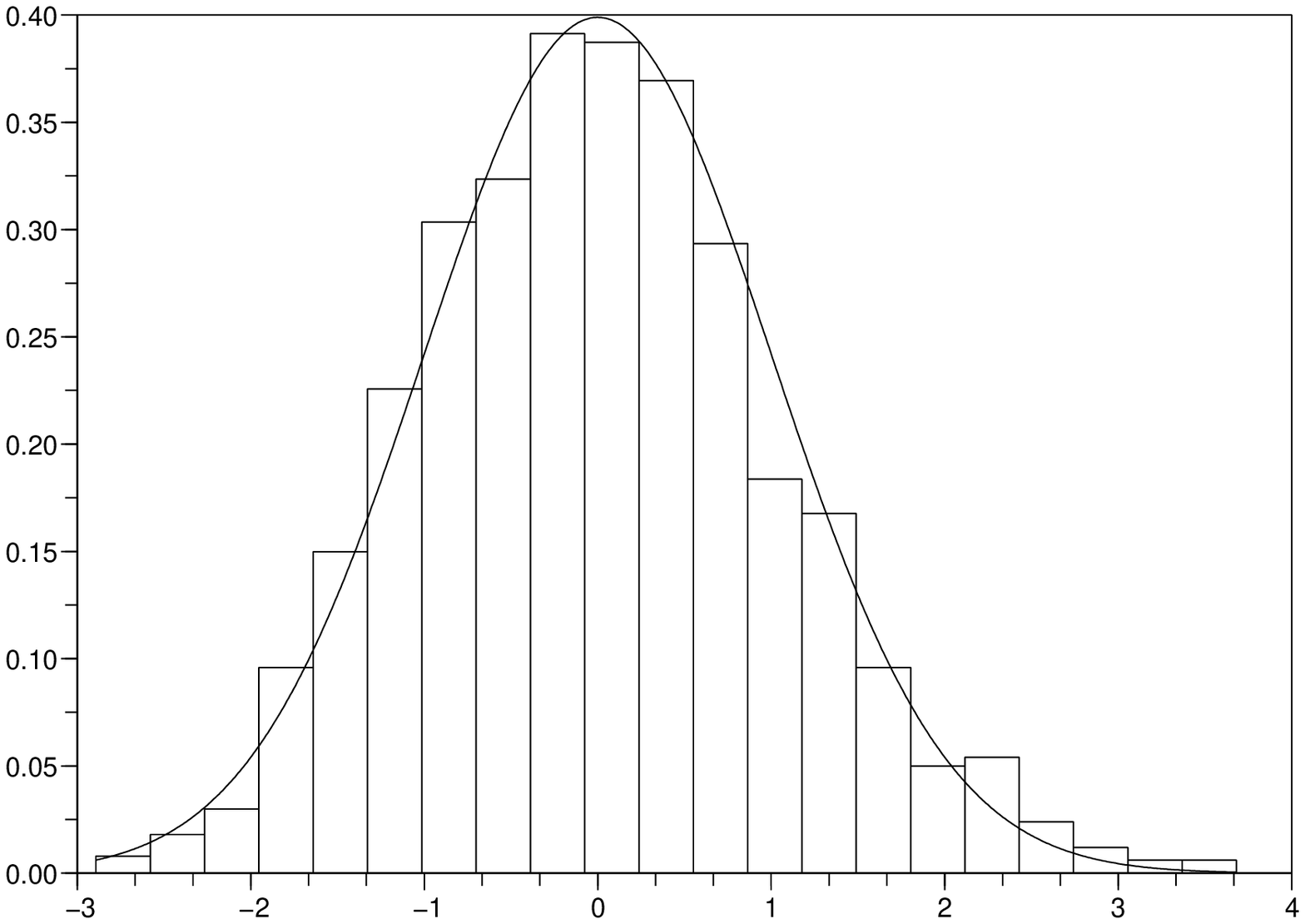}
\caption{Histogram of the kurtosis statistics over the same $1600$ simulated
fields.\label{fig-ku}}
\end{figure}

In a way of confirmation of the CLT, simulations has been performed. It has
been
chose to put $c_l=l^{-4}$. $1600$ fields were simulated on the torus with
these specifications and for each of them the skewness and kurtosis
statistics were computed for $N=2^{12}$. The values obtained were
then normalized dividing the by the theoretical standard deviations of
these statistics, computed using formulas (\ref{varb}),
(\ref{varsk}), (\ref{varkur1}) and (\ref{varkur2}).
Histograms are presented in Figures \ref{fig-sk} and \ref{fig-ku}, showing a
good accordance with Gaussianity.
\bibliography{torus}
\bibliographystyle{kluwer}

{\parindent0pt
\hbox{
\vtop{\hsize 6cm\footnotesize
PB \& DM

Dipartimento di Matematica

Universit\`a di Roma {\it Tor Vergata}

Via della Ricerca Scientifica

00161 Roma (Italy)
\medskip

+39 06 7259 4847

fax. +39 06 7259 4699

baldi@mat.uniroma2.it

marinucc@mat.uniroma2.it
}
\qquad
\vtop{\hsize 7cm
\footnotesize
GK \& DP

Laboratoire de Probabilit\'es et
Mod\`eles Al\'eatoires

2, Pl. Jussieu

75251 Paris Cedex 05, France
\medskip

+33 144277960

fax. +33 144277223

kerk@math.jussieu.fr,

picard@math.jussieu.fr
}
}
}
\end{document}